\documentclass{larson-chapter}

\begin{document}

\setcounter{chapter}{0}

\chapter{Unitary Systems, Wavelet Sets, and Operator-Theoretic Interpolation of
Wavelets and Frames}

\markboth{Larson}{Wavelet Sets and Operator-Theoretic Interpolation}

\author{David R. Larson}

\address{Department of Mathematics\\ Texas A\&M University\\
College Station, TX 77845\\
E-mail: larson@math.tamu.edu}

\begin{abstract}
A wavelet is a special case of a vector in a separable Hilbert space
that generates a basis under the action of a collection, or system,
of unitary operators. We will describe the operator-interpolation
approach to wavelet theory using the local commutant of a system.
This is really an abstract application of the theory of operator
algebras to wavelet theory. The concrete applications of this method
include results obtained using specially constructed families of
wavelet sets.   A frame is a sequence of vectors in a Hilbert space
which is a compression of a basis for a larger space.  This is not
the usual definition in the frame literature, but it is easily
equivalent to the usual definition. Because of this compression
relationship between frames and bases, the unitary system approach
to wavelets (and more generally: wandering vectors) is perfectly
adaptable to frame theory.    The use of the local commutant is
along the same lines as in the wavelet theory. Finally, we discuss
constructions of frames with special properties using targeted
decompositions of positive operators, and related problems.

\end{abstract}

\section{Introduction}

This is a write-up of of a tutorial series of three talks which I
gave as part of the "Workshop on Functional and Harmonic Analyses of
Wavelets and Frames" held August 4-7, 2004 at the National
University of Singapore. I will first give the titles and abstracts
essentially as they appeared in the workshop schedule. I will say
that the actual style of write-up of these notes will be structured
a bit differently, but only in that more than three sections will be
given, and subsections indicated, to (hopefully) improve
expositional quality.

\subsection {Talks and Abstracts}
 (a) "Unitary Systems and Wavelet Sets":  A wavelet is a special case of a
vector
 in a separable Hilbert space that generates a basis under the
 action of a collection, or "system", of unitary operators defined in terms of
 translation and dilation operations.  This approach to wavelet theory goes
 back, in particular, to earlier work of  Goodman, Lee and Tang [25] in the
 context of multiresolution analysis.  We will begin by describing the
 operator-interpolation approach to wavelet theory using the local commutant of
a
 system that was worked out by the speaker and his collaborators a few years
ago.  This
 is really an abstract application of the theory of operator algebras,
 mainly von Neumann algebras, to wavelet theory.  The concrete applications of
 operator-interpolation to wavelet theory include results obtained using
specially
 constructed families of wavelet sets.  In fact X. Dai and the speaker had
originally
 developed our theory of wavelet sets [11] specifically to take advantage of their
 natural and elegant relationships with these wavelet unitary systems. We will
also discuss
 some new results and open questions.

 (b) "Unitary Systems and Frames":  A frame is a sequence of vectors in a
Hilbert
 space which is a compression of a basis for a larger space.  (This is not the
usual
 definition in the frame literature, but it is equivalent to the usual
definition.  In
 this spirit, the usual "inequality" definition can be thought of as an abstract
 characterization of a compression of a basis.)  Because of this compression
relationship
 between frames and bases, the unitary system approach to wavelets (and more
generally:
 wandering vectors) is perfectly adaptable to frame theory.  This idea was
developed
 into a theory a few years ago by D. Han and the speaker [33].  The use of the
 local commutant is along the same lines.

(c)  "Decompositions of Operators and Operator-Valued Frames":  We
will discuss some joint work with K. Kornelson and others on
construction of frames with targeted properties [16, 42].  These are
related to targeted decompositions of positive operators.

\subsection{Some Background} It might be appropriate to give some
comments of a personal-historical nature, before continuing with the
technical aspects. My particular point of view on "wavelet theory",
which was developed jointly with my good friend and colleague Xingde
Dai, began in the summer of 1992. Before then, I was strictly an
operator theorist. I had heard my approximation theory colleagues
and friends at Texas A\&M University talk about wavelets and frames,
and Xingde had frequently mentioned these topics to me when he was
finishing his Ph.D. at A\&M (he was my student, graduating in 1990,
with a thesis [9] on the subject of \emph{nest algebras}). But it
was this meeting of minds we had in June of 1992 that was the
turning point for me. We formulated an approach to wavelet theory
(and as it turned out ultimately, to frame-wavelet theory as well)
that I felt we could "really understand" as operator algebraists.
This was the abstract \emph{unitary system} approach. Dai knew the
unitary operator approach to multiresolution analysis that had been
recently (at that time) published by Goodman, Lee and Tang 25], and
he suggested to me that we should try to go further with these ideas
in an attempt to get a some type of tractable classification of
\emph{all} wavelets. We went, in fact, in some completely different
directions. The first paper that came out of this was the AMS Memoir
[11] with Dai. The second paper was our paper [12] with Dai and
Speegle, which proved the existence of single wavelets in higher
dimensions, for arbitrary expansive dilations. After that, several
papers followed including the AMS Memoir [33] with Deguang Han, and
the Wutam Consortium paper [52], as well as the papers [10, 13, 26,
27, 28, 30, 34, 43, 44] by my students Dai, Gu, Han and Lu, and
their collaborators, and the papers [2, 39, 40] with my colleagues
Azoff, Ionascu, and Pearcy.

The paper [11] with Dai mentioned above, which was published in
1998, culminated about two years of work on this topic by the
authors. It contained our entire operator-theoretic approach to
wavelet theory, and was completed in December 1994.  This work,
while theoretical, had much hands-on experimentation in its
development, and resulted in certain theorems we were able to prove
concerning constructions of new families of wavelets. In order to
conduct successful \emph{experiments} with our operator techniques,
we needed a supply of easily computable \emph{test} wavelets: that
is, wavelets which were very amenable to \emph{paper and pencil
computations}. We discovered that certain sets, we called
\emph{wavelet sets}, existed in abundance, and we computed many
concrete examples of them along the way toward proving our results
of [11]. Several of these were given as examples in [11, Example
4.5, items $(i)\rightarrow (xi)$]. Some of these are given in
Section 2.6.1 of the present article.

Most of our work in [11] on the local commutant and the theory of
wavelet sets was accomplished in the two-month period July-September
1992. The first time Dai and I used the terms \emph{wavelet set} and
\emph{local commutant}, as well as the first time we discussed what
we referred to as the \emph{connectedness problem for wavelets}, was
in a talk in a Special Session on \emph{Operator Algebras} in the
October 1992 AMS Sectional Meeting in Dayton.

Along the way a graduate student at A\&M, Darrin Speegle, who was
enrolled in a seminar course of Larson on the manuscript of [11],
answered an open question Larson gave out in class by proving that
the set of all wavelet sets for a given wavelet system is
\emph{connected} in the symmetric difference metric on the class of
measurable sets of finite measure. That resulted in a paper [50]
which became part of his thesis (which was directed by William
Johnson of Texas A\&M), and Speegle subsequently joined forces with
Dai and Larson [12] to prove that wavelet sets (and indeed,
wavelets) exist in much greater generality than the prevailing
\emph{folklore} dictated. We received some attention for our work,
and especially we thank Guido Weiss and John Benedetto for
recognizing our work. This led to a flurry of papers by a number of
authors, notably [4, 6, 13], and also led to the paper [52] by the
Wutam Consortium, which was a group led by Guido Weiss and Larson,
consisting of 14 researchers--students and postdocs of Weiss and
Larson--based at Washington University and Texas A\&M University,
for the purpose of doing basic research on wavelet theory.

\subsubsection{Interpolation}
The main point of the \textit{operator-theoretic interpolation} of
wavelets (and frames) that Dai and I developed is that  new wavelets
can be obtained as linear combinations of known ones using
\textit{coefficients} which are not necessarily scalars but can be
taken to be  \textit{operators} (in fact, \emph{Fourier
multipliers}) in a certain class.  The ideas involved in this, and
the essential computations, all extend naturally to more general
unitary systems and \emph{wandering vectors}, and I think that much
of the theory is best-put in this abstract setting because clarity
is enhanced, and because many of the methods work for more involved
systems that are important to applied harmonic analysis, such as
Gabor and generalized Gabor systems, and various types of
\textit{frame} unitary systems.

\subsubsection{Some Basic Terminology}

This article will concern bounded linear operators on separable
Hilbert spaces. The set of all bounded linear operators on a Hilbert
space $H$ will be denoted by $B(H)$.  By a \textit{bilateral shift}
$U$ on $H$ we mean a unitary operator $U$ for which there exists a
closed linear subspace $E \subset H$ with the property that the
family of subspaces $\{U^nE: n \in \mathbb{Z}\}$ are orthogonal and
give a direct-sum decomposition of $H$.  The subspace $E$ is called
a \emph{complete wandering subspace} for $U$. The
\textit{multiplicity} of $U$ is defined to be the dimension of $E$.

The \emph{strong operator topology} on $B(H$ is the topology of
pointwise convergence, and the \emph{weak operator topology} is the
weakest topology such that the vector functionals $\omega_{x,y}$ on
$B(H)$ defined by $A \mapsto \langle Ax, y \rangle$, $A \in B(H)$,
$x,y \in H$, are all continuous.  An \emph{algebra of operators} is
a linear subspace of $B(H)$ which is closed under multiplication. An
\emph{operator algebra} is an algebra of operators which is
\emph{norm-closed}. A subset $\mathcal{S} \subset B(H)$ is called
\emph{selfadjoint} if whenever $A \in \mathcal{S}$ then also $A^*
\in \mathcal{S}$.  A $C^*$-\emph{algebra} is a self-adjoint operator
algebra. A \emph{von Neumann algebra} is a $C^*$-algebra which is
closed in the weak operator topology.  For a unital operator
algebra, it is well known that being closed in the weak operator
topology is equivalent to being closed in the closed in the strong
operator topology.

The \emph{commutant} of a set $\mathcal{S}$ of operators in $B(H)$
is the family of all operators in $B(H)$ that \emph{commute} with
every operator in $\mathcal{S}$.  It is closed under addition and
multiplication, so is an algebra.  And it is clearly closed in both
the weak operator topology and the strong operator topology.  We use
the standard \emph{prime} notation for the commutant. So the
commutant of a subset $\mathcal{S} \subset B(H)$ is denoted:~~
$\mathcal{S}^\prime := \{A \in B(H): AS = SA, ~~ S \in
\mathcal{S}\}$.

The commutant of a selfadjoint set of operators is clearly a von
Neumann algebra.  Moreover, by a famous theorem of Fuglede every
operator which commutes with a normal operator $N$ also commutes
with its adjoint $N^*$, and hence the commutant of any set of
\emph{normal} operators is also a von Neumann algebra.  So, of
particular relevance to this work, the commutant of any set of
\emph{unitary} operators is a von Neumann algebra.

One of the main tools in this work is the \emph{local commutant} of
a system of unitary operators.  (See section 2.4.)  This is a
natural generalization of the commutant of the system, and like the
commutant it is a linear space of operators which is closed in the
weak and the strong operator topologies, but unlike the commutant it
is usually not selfadjoint, and is usually not closed under
multiplication.  It contains the commutant of the system, but can be
much larger than the commutant.  The local commutant of a wavelet
unitary system captures all the information about the wavelet system
in an essential way, and this gives the \emph{flavor} of our
approach to the subject.

If $U$ is a unitary operator and $\mathcal{A}$ is an operator
algebra, then $U$ is said to \emph{normalize} $\mathcal{A}$ if
~~$U^\star \cdot \mathcal{A} \cdot U = \mathcal{A}$~.  In the most
interesting cases of operator-theoretic interpolation: that is,
those cases that yield the strongest structural results, the
relevant unitaries in the local commutant of the system normalize
the commutant of the system.

\subsubsection{Acknowledgements}
I want to take the opportunity to
thank the organizers of this wonderful workshop at the National
University of Singapore for their splendid hospitality and great
organization, and for inviting me to give the  series of
tutorial-style talks that resulted in this write-up.  I also want to
state that the work discussed in this article was supported by
grants from the United States National Science Foundation.

\section{Unitary Systems and Wavelet Sets}

We define a \textit{unitary system} to be simply a collection of
unitary operators $\mathcal U$ acting on a Hilbert space $H$ which
contains the identity operator.  The \textit{interesting} unitary
systems all have additional structural properties of various types.
  We will say that a vector $\psi \in H$ is \textit{wandering} for
$\mathcal U$ if the set
\begin{equation} \mathcal U \psi := \{U\psi : U \in \mathcal U\} \end{equation}
is an
orthonormal set, and we will call $\psi$ a \textit{complete
wandering vector} for $\mathcal U$ if $\mathcal U \psi$ spans $H$.
This (abstract) point of view can be useful.  Write $$\mathcal
W(\mathcal U)$$ for the set of complete wandering vectors for
$\mathcal U$.

\subsection{The One-Dimensional Wavelet System}
For simplicity of presentation, much of the work in this article
will deal with one-dimensional wavelets, and in particular, the
dyadic case.  The other cases: non-dyadic and in higher dimensions,
are well-described in the literature and are at least notationally
more complicated.

\subsubsection{Dyadic Wavelets} A \textit{dyadic orthonormal} wavelet in one
dimension
is a unit vector $\psi \in L^2(\mathbb{R}, \mu)$, with $\mu$
Lebesgue measure, with the property that the set
\begin{equation} \{2^{\frac{n}{2}}\psi(2^n t - l) : n,l \in
\mathbb{Z}\}\end{equation} of all integral translates of $\psi$
followed by dilations by arbitrary integral powers of $2$, is an
orthonormal basis for $L^2(\mathbb{R},\mu)$.  The term \emph{dyadic}
refers to the dilation factor "$2$".  The term \emph{mother wavelet}
is also used in the literature for $\psi$. Then the functions
$$\psi_{n,l} := 2^{\frac{n}{2}} \psi(2^n t - l)$$ are called elements
of the wavelet basis generated by the "mother".  The functions
$\psi_{n,l}$ will not themselves be mother wavelets unless $n = 0$.

Let $T$ and $D$ be the translation (by $1$) and dilation (by $2$)
unitary operators in $B(L^2(\mathbb{R})$ given by $(Tf)(t) = f(t-1)$
and $(Df)(t) = \sqrt{2}f(2t)$.  Then
$$2^{\frac{n}{2}}\psi(2^nt - l) = (D^nT^l \psi)(t)$$
for all $n,l \in \mathbb{Z}$.  Operator-theoretically, the operators
$T, D$ are \textit{bilateral shifts} of \textit{infinite
multiplicity}.  It is obvious that $L^2([0,1])$, considered as a
subspace of $L^2(\mathbb{R})$, is a complete wandering subspace for
$T$, and that  $L^2([-2, -1] \cup [1, 2])$ is a complete wandering
subspace for $D$.

\subsubsection{The Dyadic Unitary System}Let $\mathcal{U}_{D,T}$ be the unitary
system defined
by
\begin{equation} \mathcal{U}_{D,T}= \{D^nT^l : n,l \in \mathbb{Z}\}
\end{equation}
where $D$ and $T$ are the operators defined above.  Then $\psi$ is a
dyadic orthonormal wavelet if and only if $\psi$ is a complete
wandering vector for the unitary system $\mathcal{U}_{D,T}$. This
was our original motivation for developing the abstract unitary
system theory.  Write
\begin{equation} \mathcal{W}(D,T) := \mathcal{W}(\mathcal{U}_{D,T})
\end{equation}
  to denote the set of all dyadic orthonormal wavelets in one
  dimension.

An abstract interpretation is that, since $D$ is a bilateral shift
it has (many) complete wandering subspaces,  and a wavelet for the
system is a vector $\psi$ whose translation space (that is, the
closed linear span of $\{T^k: k \in \mathbb{Z}\}$ is a complete
wandering subspace for $D$.  Hence $\psi$ must generate an
orthonormal basis for the entire Hilbert space under the action of
the unitary system.

\subsubsection{Non-Dyadic Wavelets in One Dimension}  In one
dimension, there are non-dyadic orthonormal wavelets: i.e. wavelets
for all possible dilation factors besides $2$ (the dyadic case). We
said "possible", because the scales $\{0,1,-1\}$ are excluded as
scales because the dilation operators they would introduce are not
bilateral shifts. All other real numbers for scales yield wavelet
theories. In [11, Example 4.5 (x)] a family of examples is given of
three-interval wavelet sets (and hence wavelets) for all scales $d
\geq 2$, and it was noted there that such a family also exists for
dilation factors $1 < d \leq 2$. There is some recent (yet
unpublished) work that has been done, by REU students and mentors,
building on this, classifying finite-interval wavelet sets for all
possible real (positive and negative scale factors).  I mentioned
this work, in passing, in my talk.

\subsection{N dimensions}

\subsubsection{The Expansive-Dilation Case}

Let $1 \leq m < \infty$, and let $A$ be an $n \times n$ real matrix
which is \textit{expansive} (equivalently, all (complex) eigenvalues
have modulus $>1$).  By a \emph{dilation - $A$ regular-translation
orthonormal wavelet} we mean a function $\psi \in L^2(\mathbb{R}^n)$
such that
\begin{equation} \{|det(A)|^{\frac{n}{2}} \psi(A^n t - (l_1, l_2, ..., l_n)^t :
n,l \in \mathbb{Z}\} \end{equation}
where $t = (t_1, ..., t_n)^t$, is an orthonormal basis for
$L^2(\mathbb{R}^n ; m)$.  (Here $m$ is product Lebesgue measure, and
the superscript "t" means transpose.)

If $A \in M_n(\mathbb{R})$ is invertible (so in particular if $A$ is
expansive), then it is very easy to verify that the operator defined
by
\begin{equation} (D_Af)(t) = |det A|^{\frac12} f(At) \end{equation}
for $f \in L^2(\mathbb{R}^n)$, $t \in \mathbb{R}^n$, is
\emph{unitary}. For $1 \leq i \leq n$, let $T_i$ be the unitary
operator determined by translation by $1$ in the $i^{th}$ coordinate
direction.  The set (5) above is then
\begin{equation} \{D^k_A T^{l_1}_1 \cdot\cdot\cdot T^{l_n}_n \psi : k,l_i \in
\mathbb{Z}\} \end{equation}

If the dilation matrix $A$ is expansive, but the translations are
along some oblique lattice, then there is an invertible real $n
\times n$ matrix $T$ such that conjugation with $D_T$ takes the
entire wavelet system to a regular-translation expansive-dilation
matrix. This is easily worked out, and was shown in detail in [39]
in the context of working out a complete theory of unitary
equivalence of wavelet systems. Hence the wavelet theories are
equivalent.

\subsubsection{ The Non-Expansive Dilation Case}

Much work has been accomplished concerning the existance of wavelets
for dilation matrices $A$ which are not expansive.  Some of the
original work was accomplished in the Ph.D. theses of Q. Gu and  D.
Speegle, when they were together finishing up at Texas A\&M. Some
significant additional work was accomplished by Speegle in [49], and
also by others. In  [39], with Ionascu and Pearcy we proved that if
an nxn real invertible matrix $A$ is not similar (in the nxn complex
matrices) to a unitary matrix, then the corresponding dilation
operator $D_A$ is in fact a bilateral shift of infinite
multiplicity.  If a dilation matrix were to admit any type of
wavelet (or frame-wavelet) theory, then it is well-known that a
necessary condition would be that the corresponding dilation
operator would have to be a bilateral shift of infinite
multiplicity.  I am happy to report that in very recent work [45],
with E. Schulz, D. Speegle, and K. Taylor, we have succeeded in
showing that this minimal condition is in fact sufficient: such a
matrix, with regular translation lattice, admits a (perhaps
infinite) tuple of functions, which collectively generates a
frame-wavelet under the action of this unitary system.

\subsection{Abstract Systems}

\subsubsection{Restrictions on Wandering Vectors}

We note that \textit{most} unitary systems $\mathcal{U}$ do not have
complete wandering vectors.  For $\mathcal{W}(\mathcal{U})$ to be
nonempty, the set $\mathcal{U}$ must be very special. It must be
\emph{countable} if it acts separably (i.e. on a separable Hilbert
space), and it must be \emph{discrete} in the strong operator
topology because if $U,V \in \mathcal{U}$ and if $x$ is a wandering
vector for $\mathcal{U}$ then
$$\|U - V \| \geq \| Ux - Vx \| = \sqrt{2}$$
Certain other properties are forced on $\mathcal{U}$ by the presence
of a wandering vector.  One purpose of [11] was to study such
properties.  Indeed, it was a matter of some surprise to us to
discover that such a theory is viable even in some considerable
generality.  For perspective, it is useful to note that while
$\mathcal{U}_{D,T}$ has complete wandering vectors, the reversed
system $$\mathcal{U}_{T,D} = \{T^lD^n : n,l \in \mathbb{Z}\}$$
\textit{fails} to have a complete wandering vector.  (A proof of
this was given in the introduction to [11].)

\subsubsection{Group Systems}
An example which is important to the theory is the following:  let
$G$ be an arbitrary countable group, and let $H = l^2(G)$.  Let
$\pi$ be the (left) regular representation of $G$ on $H$.  Then
every element of $G$ gives a complete wandering vector for the
unitary system $$\mathcal{U} := \pi(G) .$$  (If $h \in G$ it is
clear that the vector $\lambda_h  \in  l^2(G)$, which is defined to
have $1$ in the $h$ position and $0$ elsewhere, is in
$\mathcal{W}(\mathcal{U})$.) If a unitary system is a
\textit{group}, and if it has a complete wandering vector, it is not
hard to show that it is unitarily equivalent to this example.

\subsection{The Local Commutant}

\subsubsection{ The Local Commutant of the System $\mathcal{U}_{D,T}$}

Computational aspects of operator theory can be introduced into the
wavelet framework in an elementary way.  Here is the way we
originally did it:  Fix a wavelet $\psi$ and consider the set of all
operators  $S \in  B(L^2(\mathbb{R}))$ which \textit{commute} with
the \emph{action} of dilation and translation on $\psi$.  That is,
require

\begin{equation}(S\psi)(2^nt-l) = S(\psi(2^nt - l))\end{equation}
 or equivalently
\begin{equation}D^nT^lS\psi = SD^nT^l\psi \end{equation} for all $n,l \in
\mathbb{Z}$.  Call
this the \emph{local commutant of the wavelet system
$\mathcal{U}_{D,T}$ at the vector $\psi$}. (In our first preliminary
writings and talks we called it the \emph{point commutant} of the
system . Formally, the local commutant of the dyadic wavelet system
on $L^2(\mathbb{R})$ is:
\begin{equation}\mathcal{C}_\psi(\mathcal{U}_{D,T}) := \{S \in
B(L^2(\mathbb{R})): (SD^nT^l - D^nT^lS)\psi = 0, \forall n,l \in
\mathbb{Z}\}\end{equation}
This is a linear subspace of $B(H)$ which is closed in the strong
operator topology, and in the weak operator topology, and it clearly
contains the \emph{commutant} of $\{ D,T \}$.

A motivating example is that if $\eta$ is any other wavelet, let $V
:= V_\psi^\eta$ be the unitary (we call it the \textit{interpolation
unitary}) that takes the basis $\psi_{n,l}$ to the basis
$\eta_{n,l}$.  That is, $V\psi_{n,l} = \eta_{n,l}$ for all $n,l \in
\mathbb{Z}$.  Then $\eta = V\psi$, so $VD^nT^l\psi = D^nT^lS\psi$
hence $V \in \mathcal{C}_\psi(\mathcal{U}_{D,T})$.

In the case of a pair of complete wandering vectors $\psi,\eta$ for
a general unitary system $\mathcal{U}$, we will use the same
notation $V_\psi^\eta$ for the unitary that takes the vector $U\psi$
to $U\eta$ for all $U \in \mathcal{U}$.

This simple-minded idea is reversible, so for every unitary $V$ in
$\mathcal{C}_\psi(\mathcal{U}_{D,T})$ the vector $V \psi$ is a
wavelet. This correspondence between unitaries in
$\mathcal{C}_\psi(D,T)$ and dyadic orthonormal wavelets is
one-to-one and onto (see Proposition ~1.)   This turns out to be
useful, because it leads to some new formulas relating to
decomposition and factorization results for wavelets, making use of
the \textit{linear} and \textit{multiplicative} properties of
$\mathcal{C}_\psi(D,T)$.

It turns out (a proof is required) that the entire local commutant
of the system $\mathcal{U}_{D,T}$ at a wavelet $\psi$ is \emph{not}
closed under multiplication, but it also turns out (also via a
proof) that for \emph{most} (and perhaps \emph{all}) wavelets $\psi$
the local commutant at $\psi$ contains many noncommutative operator
algebras (in fact von Neumann algebras) as subsets, and their
unitary groups \emph{parameterize} norm-arcwise-connected families
of wavelets. Moreover, $\mathcal{C}_\psi(D,T)$ is closed under
\emph{left multiplication} by the commutant $\{ D,T \}^\prime$,
which turns out to be an abelian nonatomic von Neumann algebra.  The
fact that $\mathcal{C}_\psi(D,T)$ is a \emph{left module} under $\{
D,T \}^\prime$ leads to a method of obtaining new wavelets from old,
and of obtaining connectedness results for wavelets, which we called
\textit{operator-theoretic interpolation} of wavelets in [DL], (or
simply \emph{operator-interpolation}).

\subsubsection{The Local Commutant of an Abstract Unitary System}

More generally, let $\mathcal{S} \subset B(H)$ be a set of
operators, where $H$ is a separable Hilbert space, and let $x \in H$
be a nonzero vector, and \textit{formally} define the \textit{local
commutant} of $\mathcal{S}$ at $x$ by
$$\mathcal{C}_x(\mathcal{S}) := \{A \in B(H) : (AS - SA)x = 0, S \in
\mathcal{S}\}$$

As in the wavelet case, this is a weakly and strongly closed linear
subspace of $B(H)$ which contains the commutant $\mathcal{S}^\prime$
of $\mathcal{S}$.  If $x$ is \textit{cyclic} for $\mathcal{S}$ in
the sense that span$(\mathcal{S}x)$ is dense in $H$, then $x$
\textit{separates} $\mathcal{C}_x(\mathcal{S})$ in the sense that
for $S \in \mathcal{C}_x(\mathcal{S})$, we have $Sx = 0$ iff $x =
0$.
Indeed, if $A \in \mathcal{C}_x(\mathcal{S})$ and if $Ax = 0$, then
for any $S \in \mathcal{S}$ we have $ASx = SAx = 0$, so $A
\mathcal{S}x = 0$, and hence $A = 0$.

If $A \in \mathcal{C}_x(\mathcal{S})$ and $B \in
\mathcal{S}^{\prime}$,  let $C = BA$.  Then for all $S \in
\mathcal{S}$,
$$(CS - SC)x = B(AS)x - (SB)Ax = B(SA)x - (BS)Ax = 0$$
because $ASx = SAx$ since $A \in \mathcal{C}_x(\mathcal{S})$, and
$SB = BS$ since $B \in \mathcal{S}^{\prime}$.
Hence $\mathcal{C}_x(\mathcal{S})$ is closed under left
multiplication by operators in $\mathcal{S}^{\prime}$. That is,
$\mathcal{C}_x(\mathcal{S}$ is a \emph{left module} over
$\mathcal{S}^{\prime}$.

It is interesting that, if in addition $\mathcal{S}$ is a
multiplicative semigroup, then in fact $\mathcal{C}_x(\mathcal{S})$
is identical with the commutant
$\mathcal{S}^\prime$ so in this case the commutant is not a new
structure.  To see this, suppose $A \in \mathcal{C}_x(\mathcal{S})$.
Then for each $S,T \in \mathcal{S}$ we have $ST \in \mathcal{S}$,
and so
$$AS(Tx) = (ST)Ax = S(ATx) = (S)Tx$$
So since $T \in \mathcal{S}$ was arbitrary and span$(\mathcal{S}x) =
H$, it follows that $AS = SA$.

\begin{proposition} If $\mathcal{U}$ is any unitary system for which
$\mathcal{W}(\mathcal{U}) \neq \emptyset$, then for any $\psi \in
\mathcal{W}(\mathcal{U})$
$$\mathcal{W}(\mathcal{U}) = \{U\psi : U
\textit{ is a unitary operator in }  \mathcal{C}_\psi(\mathcal{U})\}
$$ and the correspondence $U \rightarrow U\psi$ is one-to-one.
\end{proposition}

A \emph{Riesz basis} for a Hilbert space $H$ is the image under a
bounded invertible operator of an orthonormal basis.  Proposition 1
generalizes to generators of Riesz bases.  A \emph{Riesz vector} for
a unitary system $\mathcal{U}$ is defined to be a vector $\psi$ for
which $\mathcal{U}\psi := \{U\psi: U \in \mathcal{U}\}$ is a Riesz
basis for the closed linear span of $\mathcal{U}\psi$, and it is
called \emph{complete} if $\overline{\texttt{span}}$
$\mathcal{U}\psi = H$. Let $\mathcal{RW(U)}$ denote the set of all
complete Riesz vectors for $\mathcal{U}$.

\begin{proposition} Let $\mathcal{U}$ be a unitary system on a Hilbert space
$H$.  If $\psi$ is a complete Riesz vector for $\mathcal{U}$ , then
$$\mathcal{RW(U)} = \{A\psi : A \textit{ is an operator in }
\mathcal{C}_\psi(\mathcal{U}) \textit{ that is invertible in }
B(H)\}.$$
\end{proposition}

\subsubsection{Operator-Theoretic Interpolation }

Now suppose $\mathcal{U}$ is a unitary system, such as
$\mathcal{U}_{D,T}$, and suppose $\{\psi_1, \psi_2, \dots, \psi_m \}
\subset \mathcal{W}(\mathcal{U})$.  (In the case of
$\mathcal{U}_{D,T}$, this means that $(\psi_1, \psi_2, \dots, \psi_n
)$ is an n-tuple of wavelets.

Let $(A_1, A_2, \dots, A_n)$ be an n-tuple of operators in the
commutant  $\mathcal{U}^{\prime}$ of $\mathcal{U}$, and let $\eta$
be the vector
$$ \eta := A_1 \psi_1 + A_2 \psi_2 + \dots + A_n \psi_n ~.$$
Then

$$ \eta = A_1 \psi_1 + A_2 V_{\psi_1}^{\psi_2}\psi_1 + \dots A_n
V_{\psi_1}^{\psi_n}\psi_1 $$
\begin{equation} ~~= (A_1 + A_2
V_{\psi_1}^{\psi_2} + \dots + A_n V_{\psi_1}^{\psi_n}) \psi_1~~.
\end{equation}

We say that $\eta$ is obtained by $\emph{operator interpolation}$
from $\{\psi_1, \psi_2, \dots, \psi_m \}$.  Since
$\mathcal{C}_{\psi_1}(\mathcal{U})$ is a left $\mathcal{U}^{\prime}$
- module, it follows that the operator
\begin{equation} A := A_1 + A_2 V_{\psi_1}^{\psi_2} + \dots A_n
V_{\psi_1}^{\psi_n}
\end{equation}
is an element of $\mathcal{C}_{\psi_1}(\mathcal{U})$.  Moreover, if
$B$ is another element of $\mathcal{C}_{\psi_1}(\mathcal{U})$ such
that $\eta = B \psi_1$, then $~A - B ~\in
\mathcal{C}_{\psi_1}(\mathcal{U})$ and $(A - B)\psi_1 = 0$.  So
since $\psi_1$ \em{separates} $\mathcal{C}_{\psi_1}(\mathcal{U})$ it
follows that $A = B$.  Thus $A$ is the \em{unique} element of
$\mathcal{C}_{\psi_1}(\mathcal{U})$ that takes $\psi_1$ to $\eta$.
Let $\mathcal{S}_{\psi_1, \dots, \psi_n}$ be the family of all
finite sums of the form
$$ \sum^n_{i=0} A_i V_{\psi_1}^{\psi_i} ~~~.$$
The is the left module of $\mathcal{U}^{\prime}$ generated by $\{I,
V_{\psi_1}^{\psi_2}, \dots, V_{\psi_1}^{\psi_n}\}$.  It is the
$\mathcal{U}^{\prime}$-\emph{linear span} of $\{I,
V_{\psi_1}^{\psi_2}, \dots, V_{\psi_1}^{\psi_n}\}$. Let
\begin{equation}
\mathcal{M}_{\psi_1, \dots, \psi_n} := (\mathcal{S}_{\psi_1, \dots,
\psi_n}) \psi_1
\end{equation}

\noindent So
\[ \mathcal{M}_{\psi_1, \dots, \psi_n} ~=~ \left\{\sum^n_{i=0} A_i
\psi_i ~: ~A_i \in \mathcal{U}^{\prime}\right\} ~~~.\]

\noindent We call this the \emph{interpolation space} for
$\mathcal{U}$ generated by $(\psi_1, \dots, \psi_n)$.  From the
above discussion, it follows that for every vector $\eta \in
\mathcal{M}_{\psi_1, \psi_2, \dots, \psi_n}$ there exists a unique
operator $A \in \mathcal{C}_{\psi_1}(\mathcal{U})$ such that $\eta =
A \psi_1$, and moreover this $A$ is an element of
$\mathcal{S}_{\psi_1, \dots, \psi_n}$.

\subsubsection{Normalizing the Commutant}
In certain essential cases (and we are not sure how general this
type of case is) one can prove that an interpolation unitary
$V_\psi^\eta$ \emph{normalizes} the commutant $\mathcal{U}^{\prime}$
of the system in the sense that $V_\eta^\psi \mathcal{U}^{\prime}
V_\psi^\eta = \mathcal{U}^{\prime}$.  (Here, it is easily seen that
$(V_\psi^\eta)^* = V_\eta^\psi$.)  Write $V := V_\psi^\eta$.  If $V$
normalizes $\mathcal{U}^\prime$, then the algebra, before norm
closure, generated by
$\mathcal{U}^\prime$ and $V$ is the set of all finite sums (trig
polynomials) of the form $\sum A_n V^n$, with coefficients $A_n \in
\ \mathcal{U}^\prime$, $n \in \mathbb{Z}$. The closure in the strong
operator topology is a von Neumann algebra. Now suppose further that
{\em every power\/} of $V$ is contained in
$\mathcal{C}_\psi(\mathcal{U})$. This occurs only in special cases,
yet it occurs frequently enough to yield some general methods. Then
since $\mathcal{C}_\psi(\mathcal{U})$ is a SOT-closed linear
subspace which is closed under left multiplication by
$\mathcal{U}^\prime$, this von Neumann algebra is contained in
$\mathcal{C}_\psi(\mathcal{U})$, so its unitary group parameterizes
a norm-path-connected subset of $\mathcal{W}(\mathcal{U})$ that
contains $\psi$ and $\eta$ via the correspondence $U\to U\psi$.

In the special case of \emph{wavelets}, this is the basis for the
work that Dai and I did in [11, Chapter 5] on operator-theoretic
interpolation of wavelets.  In fact, we specialized there and
\emph{reserved} the term \emph{operator-theoretic interpolation} to
refer explicitely to the case when the interpolation unitaries
normalize the commutant. In some subsequent work, we \emph{loosened}
this restriction yielding our more general definition given in this
article, because there are cases of interest in which we weren't
able to prove normalization.  However, it turns out that if $\psi$
and $\eta$ are $s$-elementary wavelets (see section 2.5.4), then
indeed $V^\eta_\psi$ normalizes $\{D,T\}'$. (See Proposition 14.)
Moreover, $V^\eta_\psi$ has a very special form:\ after conjugating
with the Fourier transform, it is a composition operator with symbol
$a$ natural and very computable measure-preserving transformation of
$\mathbb{R}$. In fact, it is precisely this special form for
$V^\eta_\psi$ that allows us to make the computation that it
normalizes $\{D,T\}'$. On the other hand, we know of no pair
$(\psi,\eta)$ of wavelets for which $V^\eta_\psi$ fails to normalize
$\{D,T\}'$. The difficulty is simply that in general it is very hard
to do the computations.
\medskip

{\bf Problem:} If $\{\psi,\eta\}$ is a pair of dyadic orthonormal
wavelets, does the interpolation unitary $V^\eta_\psi$ normalize
$\{D,T\}'$? As mentioned above, the answer is yes if $\psi$ and
$\eta$ are $s$-elementary wavelets.

\subsubsection{An Elementary Interpolation Result}

The following result is the most elementary case of
operator-theoretic interpolation.

\begin{proposition} Let $\mathcal{U}$ be a unitary system on a
Hilbert space $H$.  If $\psi_1$ and $\psi_2$ are in
$\mathcal{W(U)}$, then $$\psi_1 + \lambda \psi_2 \in
\mathcal{RW(U)}$$ for all complex scalars $\lambda$ with $|\lambda|
\neq 1$.  More generally, if $\psi_1$ and $\psi_2$ are in
$\mathcal{RW(U})$ then there are positive constants \space $b > a
 > 0$ such that $\psi_1 + \lambda \psi_2 \in
\mathcal{RW(U)}$ for all $\lambda \in \mathbb{C}$ with either
$|\lambda| < a$ or with $|\lambda| > b$.

\end{proposition}

\begin{proof}  If $\psi_1, \psi_2 \in \mathcal{W(U)}$, let $V$ be
the unique unitary in $\mathcal{C}_{\psi_2}(\mathcal{U})$ given by
Proposition 1 such that $V\psi_2 = \psi_1$.  Then $$\psi_1 + \lambda
\psi_2 = (V + \lambda I)\psi_2 .$$  Since $V$ is unitary, $(V +
\lambda I)$ is an invertible element of
$\mathcal{C}_{\psi_2}(\mathcal{U})$ if $|\lambda | \neq 1$, so the
first conclusion follows from Proposition 2.  Now assume $\psi_1,
\psi_2 \in \mathcal{RW(U)}$.  Let $A$ be the unique invertible
element of $\mathcal{C}_{\psi_2}(\mathcal{U})$ such that $A \psi_2 =
\psi_1$, and write $\psi_1 + \lambda \psi_2 = (A + \lambda
I)\psi_2$.  Since $A$ is bounded and invertible there are $b > a
> 0$ such that $$\sigma (A) \subseteq \{z \in \mathbb{C}: a < |z| < b
\}$$ where $\sigma (A)$ denotes the spectrum of $A$ , and the same
argument applies.

\end{proof}

\subsubsection{Interpolation Pairs of Wandering Vectors}

In some cases where a pair $\psi, \eta$ of vectors in
$\mathcal{W(U)}$ are given  it turns out that the unitary $V$ in
$\mathcal{C}_\psi(\mathcal{U})$ with $V\psi = \eta$ happens to be a
\emph{symmetry} (i.e. $V^2 = I$).  Such pairs are called
\emph{interpolation pairs} of wandering vectors, and in the case
where $\mathcal{U}$ is a wavelet system, they are called
interpolation pairs of wavelets. Interpolation pairs are more
prevalent in the theory, and in particular the wavelet theory, than
one might expect. In this case (and in more complex generalizations
of this) certain linear combinations of complete wandering vectors
are themselves complete wandering vectors -- not simply complete
Riesz vectors.

\begin{proposition} Let $\mathcal{U}$ be a unitary system, let
$\psi, \eta \in \mathcal{W(U)}$, and let $V$ be the unique operator
in $\mathcal{C}_\psi(\mathcal{U})$ with $V \psi = \eta$.  Suppose
$$V^2 = I.$$ Then
$$ \cos \alpha \cdot \psi \textit{ +  } i \sin \alpha \cdot \eta \in
\mathcal{W(U)}$$ for all $0 \leq \alpha \leq 2\pi$.

\end{proposition}

The above result can be thought of as the \emph{prototype} of our
operator-theoretic interpolation results.  It is the second most
elementary case.  More generally, the scalar $\alpha$ in Proposition
4 can be replaced with an appropriate \emph{self-adjoint operator}
in the commutant of $\mathcal{U}$.  In the wavelet case, after
conjugating with the Fourier transform, which is a unitary operator,
this means that $\alpha$ can be replaced with a wide class of
nonnegative dilation-periodic (see definition below) bounded
measurable functions on $\mathbb{R}$.

\subsubsection{A Test For Interpolation Pairs}
\medskip
The following converse to Proposition 4 is typical of the type of
computations encountered in some wandering vector proofs.

\begin{proposition} Let $\mathcal{U}$ be a unitary system, let
$\psi, \eta \in \mathcal{W(U)}$, and let $V$ be the unique unitary
in $\mathcal{C}_\psi(\mathcal{U})$ with $V \psi = \eta$.  Suppose
for some $0 < \alpha < \frac{\pi}{2}$ the vector $$\rho := \cos
\alpha \cdot \psi  \textit{ +  } i \sin \alpha \cdot \eta $$ is
contained in $\mathcal{W(U)}$.  Then $$V^2 = I .$$
\end{proposition}

\begin{proof}
Since $\mathcal{U}\psi$ is a basis it will be enough to show that
$VU_1 \psi = V^{\star} U_1 \psi$ for all $U_1 \in \mathcal{U}$. So
it will suffice to prove that for all $U_1, U_2 \in \mathcal{U}$ we
have
$$\langle VU_1 \psi, U_2 \psi \rangle = \langle V^{\star}U_1 \psi,
U_2 \psi \rangle.$$
Using the fact that $V$ locally commutes with $\mathcal{U}$ at
$\psi$ we have
$$ \langle V U_1 \psi , U_2 \psi \rangle = \langle U_1
V \psi, U_2 \psi
\rangle = \langle U_1 \eta, U_2 \psi \rangle \verb"   "
\texttt{and}$$
$$ \langle V^{\star}U_1 \psi , U_2 \psi \rangle = \langle U_1
 \psi, VU_2 \psi \rangle = \langle U_1 \psi, U_2 V \psi
 \rangle = \langle U_1 \psi, U_2 \eta \rangle .$$

So we must show that $\langle U_1 \eta, U_2 \psi \rangle = \langle
U_1 \psi, U_2 \eta \rangle$  for all $U_1, U_2 \in \mathcal{U}$.

Write $\rho := \rho_{\alpha}$.  By hypothesis $\psi, \eta$ and
$\rho$ are unit vectors.  So compute
$$ 1 = \langle \rho, \rho \rangle = {\cos}^2 \alpha \cdot \langle
\psi , \psi \rangle \textit{  +  } i \sin \alpha \cos \alpha \cdot
\langle  \eta , \psi \rangle $$
$$\verb"             -  " i  \sin \alpha \cos
\alpha \cdot  \langle \psi , \eta \rangle   \textit{  +  }  {\sin}^2
\alpha \cdot \langle \eta , \eta \rangle $$
$$ \verb"        " = 1 \verb" + " i \sin \alpha \cos \alpha \cdot (\langle \eta
,
\psi \rangle \verb" - " \langle \psi , \eta \rangle ).$$

Thus, since $\sin \alpha \cos \alpha \neq 0$, we must have $\langle
\eta , \psi \rangle = \langle \psi , \eta \rangle$.  Also, for $U_1
, U_2  \in \mathcal{U}$ with $U_1 \neq U_2$ we have
$$ 0 = \langle U_1 \rho , U_2 \rho \rangle = {\cos}^2 \alpha \cdot
\langle U_1 \psi , U_2 \psi \rangle \texttt{  +  } i \sin \alpha
\cos \alpha \cdot \langle U_1 \eta , U_2 \psi \rangle $$
$$ \texttt{ - } i \sin \alpha \cos \alpha \cdot \langle U_1 \psi ,
U_2 \eta \rangle  \texttt{  +  }  {\sin}^2 \alpha \cdot \langle U_1
\eta , U_2 \eta \rangle $$
$$ = i \sin \alpha \cos \alpha \cdot ( \langle U_1 \eta , U_2 \psi
\rangle  \texttt{  -  } \langle U_1 \psi , U_2 \eta \rangle ), $$
which implies $\langle U_1 \eta , U_2 \psi \rangle = \langle U_1
\psi , U_2 \eta \rangle$ as required.

\end{proof}

The above result gives an \emph{experimental method} of checking
whether $V^2 = I$ for a given pair $\psi, \eta \in \mathcal{W(U)}$.
One just checks whether $$\rho := \frac{1}{\sqrt{2}} \psi \texttt{ +
}  \frac{i}{\sqrt{2}} \eta $$  is an element of   $\mathcal{W(U)}$,
which is much simpler than attempting to work with the infinite
matrix of $V$ with respect to the basis $\mathcal{U} \psi $ (or some
other basis for $H$).

\subsubsection{Connectedness}
If we consider again the example of the left regular representation
$\pi$ of a group $G$ on $H := l^2(G)$, then the local commutant of
$\mathcal{U} := \pi (G)$ at a vector $\psi \in \mathcal{W}(\pi (G))$
is just the commutant of $\pi (G)$. So since the unitary group of
the von Neumann algebra $(\pi (G))^{\prime}$ is
norm-arcwise-connected, it follows that $\mathcal{W} (\pi (G))$ is
norm-arcwise-connected.

Problem A in [11] asked whether $\mathcal{W}(D, L)$ is
norm-arcwise-connected.  It turned out that this conjecture was also
formulated independently by Guido Weiss ([38], [37]) from a harmonic
analysis point of view (our point of view was purely functional
analysis), and this problem (and related problems) was the primary
stimulation for the creation for the creation of the WUTAM
CONSORTIUM -- a team of 14 researchers based at Washington
University and Texas A\&M University.  (See [52].)

This \emph{connectedness conjecture} was answered {yes} in [52] for
the special case of the family of dyadic orthonormal MRA wavelets in
$L^2(\mathbb{R})$, but still remains open for the family of
\emph{arbitrary} dyadic orthonormal wavelets in $L^2(\mathbb{R})$.

In the wavelet case $\mathcal{U}_{D,T}$ , if
$\psi \in \mathcal{W}(D,T)$ then it turns out that
$\mathcal{C}_{\psi}(\mathcal {U}_{D,T})$ is in fact $\emph{much
larger}$ than $(\mathcal{U}_{D,T})^{\prime} = \{D,T\}^{\prime}$ ,
underscoring the fact that $\mathcal{U}_{D,T}$ is NOT a group.  In
particular, $\{D,T \}^{\prime}$ is abelian while $\mathcal{C}_ \psi
(\mathcal{D,T})$ is nonabelian for every wavelet $\psi$.  (The proof
of these facts are contained in [11].)

\subsection{Wavelet Sets}

\subsubsection{The Fourier Transform}

We will use the following form of the Fourier--Plancherel transform
$\mathcal{F}$ on $\mathcal{H} = L^2(\mathbb{R})$, because it is a
form \emph{normalized} so it is a unitary transformation. Although
there is another such \emph{normalized} form that is frequently
used, and actually simpler, the present form is the one we used in
our original first paper [11] involving operator theory and
wavelets, and so we will stick with it in these notes to avoid any
confusion to a reader of both.

If $f,g\in L^1(\mathbb{R}) \cap L^2(\mathbb{R})$ then
\begin{equation}\label{eq18}
(\mathcal{F}f)(s) := \frac1{\sqrt{2\pi}} \int_{\mathbb{R}} e^{-ist} f(t)dt :=
\hat f(s),
\end{equation}
and
\begin{equation}\label{eq19}
(\mathcal{F}^{-1}g)(t) = \frac1{\sqrt{2\pi}} \int_{\mathbb{R}} e^{ist}g(s)ds.
\end{equation}
We have
\[
(\mathcal{F} T_\alpha f)(s) = \frac1{\sqrt{2\pi}} \int_{\mathbb{R}} e^{-ist}
f(t-\alpha)dt = e^{-is\alpha} (\mathcal{F}f)(s).
\]
So $\mathcal{F} T_\alpha \mathcal{F}^{-1} g=e^{-is\alpha}g$. For
$A\in\mathcal{B}(\mathcal{H})$ let $\hat A$ denote
$\mathcal{F}A\mathcal{F}^{-1}$. Thus
\begin{equation}\label{eq20}
\widehat T_\alpha = M_{e^{-i\alpha s}},
\end{equation}
where for $h\in L^\infty$ we use $M_h$ to denote the multiplication operator
$f\to hf$. Since $\{M_{e^{-i\alpha s}}\colon \ \alpha\in\mathbb{R}\}$ generates
the m.a.s.a.\ $\mathcal{D}(\mathbb{R}) := \{M_h\colon \ h\in
L^\infty(\mathbb{R})\}$ as a von Neumann algebra, we have
\[
\mathcal{F}\mathcal{A}_T \mathcal{F}^{-1} = \mathcal{D}(\mathbb{R}).
\]
Similarly,
\begin{align*}
(\mathcal{F}D^nf)(s) &= \frac1{\sqrt{2\pi}} \int_{\mathbb{R}} e^{-ist} (\sqrt
2)^n f(2^nt)dt\\
&= (\sqrt 2)^{-n}\cdot \frac1{\sqrt{2\pi}} \int_{\mathbb{R}} e^{-i2^{-n}st}
f(t)dt\\
&= (\sqrt 2)^{-2} (\mathcal{F}f)(2^{2^{-n}s}) = (D^{-n}\mathcal{F}f)(s).
\end{align*}
So $\widehat D^n = D^{-n} = D^{*n}$. Therefore,
\begin{equation}\label{eq21}
\widehat D = D^{-1} = D^*.
\end{equation}

\subsubsection{The Commutant of $\{D,T\}$}
We have $\mathcal{F}\{D,T\}'\mathcal{F}^{-1} = \{\widehat D,\widehat
T\}'$. It turns out that $\{\widehat D,\widehat T\}'$ has an easy
characterization.

\begin{theorem}\label{thm6}
\[
\{\widehat D,\widehat T\}' = \{M_h\colon \ h\in L^\infty(\mathbb{R}) \text{ and
} h(s) = h(2s) \text{ a.e.}\}.
\]
\end{theorem}

\begin{proof}
Since $\widehat D = D^*$ and $D$ is unitary, it is clear that $M_h \in
\{\widehat D, \widehat T\}'$ if and only if $M_h$ commutes with $D$. So let
$g\in L^2(\mathbb{R})$ be arbitrary. Then (a.e.) we have
\begin{align*}
(M_hDg)(s) &= h(s) (\sqrt 2\ g(2s)),\quad \text{and}\\
(DM_hg)(s) &= D(h(s) g(s)) = \sqrt h(2s) g(2s).
\end{align*}
Since these must be equal a.e.\ for arbitrary $g$, we must have $h(s) = h(2s)$
a.e.
\end{proof}

Now let $E = [-2,-1)\cup [1,2)$, and for $n\in \mathbb{Z}$ let $E_n=
\{2^nx\colon \ x\in E\}$. Observe that the sets $E_n$ are disjoint and have
union $\mathbb{R}\backslash\{0\}$. So if $g$ is any uniformly bounded function
on
$E$, then $g$ extends uniquely (a.e.) to a function $\tilde g\in
L^\infty(\mathbb{R})$ satisfying
\[
\tilde g(s) = \tilde g(2s),\qquad s\in \mathbb{R},
\]
by setting
\[
\tilde g(2^ns) = g(s),\qquad s\in E, n\in \mathbb{Z},
\]
and $\tilde g(0)=0$. We have $\|\tilde g\|_\infty = \|g\|_\infty$. Conversely,
if $h$ is any  function satisfying $h(s) = h(2s)$ a.e., then $h$ is uniquely
(a.e.) determined by its restriction to $E$. This 1-1 mapping $g\to M_{\tilde
g}$ from $L^\infty(E)$ onto $\{\widehat D, \widehat T\}'$ is a $*$-isomorphism.

We will refer to a function $h$ satisfying $h(s) = h(2s)$ a.e.\ as a
2-{\em dilation periodic function}. This gives a simple algorithm
for computing a large class of wavelets from a given one, by simply
modifying the \emph{phase}:
\begin{align}
&\text{\em Given $\psi$,, let $\widehat \psi = \mathcal{F}(\psi)$, choose a
real-valued function } h\in L^\infty(E)\nonumber\\
\label{eq22}
&\text{\em arbitrarily, let $g = \exp(ih)$, extend to a 2-dilation periodic}\\
&\text{\em function $\tilde g$ as above, and compute } \psi_{\tilde g} =
\mathcal{F}^{-1}(\tilde g\widehat\psi).\nonumber
\end{align}

In the description above, the set $E$ could clearly be replaced with
$[-2\pi,-\pi)\cup [\pi,2\pi)$, or with any other ``dyadic'' set $[-2a,a)\cup
[a,2a)$ for some $a>0$.

\subsubsection{Wavelets of Computationally Elementary Form}

We now give an account of $s$-elementary and $MSF$-wavelets. The two most
elementary dyadic orthonormal wavelets are the {\em Haar wavelet\/} and {\em
Shannon's wavelet\/} (also called the Littlewood--Paley wavelet).

The Haar wavelet is the function
\begin{equation}\label{eq32}
\psi_H(t) = \left\{\begin{array}{rl}
1,&0\le t<\frac12\\
-1,&\frac12\le t\le 1\\
0,&\text{otherwise.}\end{array}\right.
\end{equation}
In this case it is very easy to see that the dilates/translates
\[
\{2^{\frac{n}2} \psi_H(2^n-\ell)\colon \ n,\ell\in \mathbb{Z}\}
\]
are orthonormal, and an elementary argument shows that their span is dense in
$L^2(\mathbb{R})$.

Shannon's wavelet is the $L^2(\mathbb{R})$-function with Fourier transform
$\widehat\psi_S = \frac1{\sqrt{2\pi}} \chi_{E_0}$ where
\begin{equation}\label{eq33}
E_0 = [-2\pi, -\pi) \cup [\pi,2\pi).
\end{equation}
The argument that $\widehat\psi_S$ is a wavelet is in a way even more
transparent than for the Haar wavelet. And it has the advantage of generalizing
nicely. For a simple argument, start from the fact that the exponents
\[
\{e^{i\ell s}\colon \ n\in \mathbb{Z}\}
\]
restricted to $[0,2\pi]$ and normalized by $\frac1{\sqrt{2\pi}}$ is an
orthonormal basis for $L^2[0,2\pi]$. Write $E_0 = E_-\cup E_+$ where $E_- =
[-2\pi, -\pi)$, $E_+ = [\pi,2\pi)$. Since $\{E_- +2\pi, E_+\}$ is a partition of
$[0,2\pi)$ and since the exponentials $e^{i\ell s}$ are invariant under
translation by $2\pi$, it follows that
\begin{equation}\label{eq34}
\left\{\frac{e^{i\ell s}}{\sqrt{2\pi}}\Big|_{E_0}\colon \ n\in
\mathbb{Z}\right\}
\end{equation}
is an orthonormal basis for $L^2(E_0)$. Since $\widehat T = M_{e^{-is}}$, this
set can be written
\begin{equation}\label{eq35}
\{\widehat T^\ell \widehat\psi_s\colon \ \ell\in \mathbb{Z}\}.
\end{equation}
Next, note that any ``dyadic interval'' of the form $J = [b,2b)$, for some $b>0$
has the property that $\{2^nJ\colon \ n\in\mathbb{Z}\}$, is a partition of
$(0,\infty)$. Similarly, any set of the form
\begin{equation}\label{eq36}
\mathcal{K} = [-2a,-a)\cup [b,2b)
\end{equation}
for $a,b>0$, has the property that
\[
\{2^n\mathcal{K}\colon \ n\in \mathbb{Z}\}
\]
is a partition of $\mathbb{R}\backslash\{0\}$. It follows that the space
$L^2(\mathcal{K})$, considered as a subspace of $L^2(\mathbb{R})$, is a complete
wandering subspace for the dilation unitary $(Df)(s) = \sqrt 2\ f(2s)$. For each
$n\in \mathbb{Z}$,
\begin{equation}\label{eq37}
D^n(L^2(\mathcal{K})) = L^2(2^{-n}\mathcal{K}).
\end{equation}
So $\bigoplus_n D^n(L^2(\mathcal{K}))$ is a direct sum decomposition of
$L^2(\mathbb{R})$. In particular $E_0$ has this property. So
\begin{equation}\label{eq38}
D^n\left\{\frac{e^{i\ell s}}{\sqrt{2\pi}}\Big|_{E_0}\colon \ \ell\in
\mathbb{Z}\right\} = \left\{\frac{e^{2^ni\ell s}}{\sqrt{2\pi}}\Big|_{2^{-n}E_0}
\colon \ \ell\in \mathbb{Z}\right\}
\end{equation}
is an orthonormal basis for $L^2(2^{-n}E_0)$ for each $n$. It follows that
\[
\{D^n\widehat T^\ell \widehat\psi_s\colon \ n,\ell\in \mathbb{Z}\}
\]
is an orthonormal basis for $L^2(\mathbb{R})$. Hence $\{D^nT^\ell\psi_s\colon \
n,\ell\in \mathbb{Z}\}$ is an orthonormal basis for $L^2(\mathbb{R})$, as
required.

The Haar wavelet can be generalized, and in fact Daubechie's
well-known continuous compactly-supported wavelet is a
generalization of the Haar wavelet. However, known generalization of
the Haar wavelet are all more complicated and difficult to work with
in hand-computations.

For our work, in order to proceed with developing an operator
algebraic theory that had a chance of directly impacting concrete
function-theoretic wavelet theory we needed a large supply of
examples of wavelets which were elementary enough to work with.
First, we found another ``Shannon-type'' wavelet in the literature.
This was the Journe wavelet, which we found described on p.~136 in
Daubechies book [14]. Its Fourier transform is $\widehat \psi_J =
\frac1{\sqrt{2\pi}} \chi_{E_J}$, where
\[
E_J = \left[-\frac{32\pi}7, -4\pi\right) \cup \left[-\pi,
-\frac{4\pi}7\right)\cup \left[\frac{4\pi}7, \pi\right) \cup \left[4\pi,
\frac{32\pi}7\right).
\]
Then, thinking the old adage ``where there's smoke there's fire!'', we
painstakingly worked out many more examples. So far, these are the basic
building
blocks in the {\em concrete\/} part of our theory. By this we mean the part of
our theory that has had some type of direct impact on function-theoretic wavelet
theory.

\subsubsection{Definition of Wavelet Set}
We define a {\em wavelet set\/} to be a measurable subset $E$ of
$\mathbb{R}$ for which $\frac1{\sqrt{2\pi}} \chi_E$ is the Fourier
transform of a wavelet. The wavelet $\widehat\psi_E :=
\frac1{\sqrt{2\pi}}\chi_E$ is called $s$-{\em elementary\/} in [11].

It turns out that this class of wavelets was also discovered and
systematically explored completely independently, and in about the
same time period, by Guido Weiss (Washington University), his
colleague and former student E. Hernandez (U.\ Madrid), and his
students X.\ Fang and X.\ Wang. In [17,37, 38] they are called MSF
(minimally supported frequency) wavelets. In signal processing, the
parameter $s$, which is the independent variable for $\widehat\psi$,
is the {\em frequency\/} variable, and the variable $t$, which is
the independent variable for $\psi$, is the {\em time\/} variable.
No function with support a subset of a wavelet set $E$ of strictly
smaller measure can be the Fourier transform of a wavelet.

\noindent {\bf Problem.} Must the support of the Fourier transform of a
wavelet contain a wavelet set? This question is open for dimension 1. It makes
sense for any finite dimension.

\subsubsection{The Spectral Set Condition}


>From the argument above describing why Shannon's wavelet is, indeed, a wavelet,
it is clear that {\em sufficient\/} conditions for $E$ to be a wavelet set are

\begin{quote}
(i)~~the normalized exponential $\frac1{\sqrt{2\pi}} e^{i\ell s}$, $\ell\in
\mathbb{Z}$, when restricted to $E$ should constitute an orthonormal basis for
$L^2(E)$ (in other words $E$ is a {\em spectral set\/} for the integer lattice
$\mathbb{Z}$),
\end{quote}

\noindent and

\begin{quote}
(ii)~~The family $\{2^nE\colon \ n\in\mathbb{Z}\}$ of dilates of $E$ by integral
powers of 2 should constitute a measurable partition (i.e.\ a partition modulo
null sets) of $\mathbb{R}$.
\end{quote}

\noindent These conditions are also necessary. In fact if a set $E$ satisfies
(i), then for it to be a wavelet set it is obvious that (ii) must be satisfied.
To show that (i) must be satisfied by a wavelet set $E$, consider the vectors
\[
\widehat D^n \widehat\psi_E = \frac1{\sqrt{2\pi}} \chi_{2^{-n}E},\qquad n\in
\mathbb{Z}.
\]
Since $\widehat\psi_E$ is a wavelet these must be orthogonal, and so the sets
$\{2^nE\colon \ n\in~\mathbb{Z}\}$ must be disjoint modulo null sets. It follows
that $\{\frac1{\sqrt{2\pi}} e^{i\ell s}|_E\colon \ \ell\in \mathbb{Z}\}$ is not
only an orthonormal set of vectors in $L^2(E)$, it must also {\em span\/}
$L^2(E)$.

It is known from the theory of \emph{spectral sets} (as an
elementary special case) that a measurable set $E$ satisfies (i) if
and only if it is a generator of a measurable partition of
$\mathbb{R}$ under translation by $2\pi$ (i.e.\ iff $\{E+2\pi
n\colon \ n\in \mathbb{Z}\}$ is a measurable partition of
$\mathbb{R}$). This result generalizes to spectral sets for the
integral lattice in $\mathbb{R}^n$. For this elementary special case
a direct proof is not hard.

\subsubsection{Translation and Dilation Congruence}

We say that measurable sets $E,F$ are {\em translation congruent modulo\/}
$2\pi$ if there is a measurable bijection $\phi\colon \  E\to F$ such that
$\phi(s)-s$ is an integral multiple of $2\pi$ for each $s\in E$; or
equivalently, if there is a measurable partition $\{E_n\colon \ n\in
\mathbb{Z}\}$ of $E$ such that
\begin{equation}\label{eq39}
\{E_n  + 2n\pi\colon \ n\in \mathbb{Z}\}
\end{equation}
is a measurable partition of $F$. Analogously, define measurable sets $G$ and
$H$
to be {\em dilation congruent modulo\/} 2 if there is a measurable bijection
$\tau\colon \ G\to H$ such that for each $s\in G$ there is an integer $n$,
depending on $s$, such that $\tau(s) = 2^ns$; or equivalently, if there is a
measurable partition $\{G_n\}^\infty_{-\infty}$ of $G$ such that
\begin{equation}\label{eq40}
\{2^nG\}^\infty_{-\infty}
\end{equation}
is a measurable partition of $H$. (Translation and dilation congruency modulo
other positive numbers of course make sense as well.)

The following lemma is useful.

\begin{lemma}\label{lem7}
Let $f\in L^2(\mathbb{R})$, and let $E = \text{\rm supp}(f)$. Then $f$ has the
property that
\[
\{e^{ins}f\colon \ n\in \mathbb{Z}\}
\]
is an orthonormal basis for $L^2(E)$ if and only if
\begin{itemize}
\item[(i)] $E$ is congruent to $[0,2\pi)$ modulo $2\pi$, and
\item[(ii)] $|f(s)| = \frac1{\sqrt{2\pi}}$ a.e.\ on $E$.
\end{itemize}
\end{lemma}

If $E$ is a measurable set which is $2\pi$-translation congruent to $[0,2\pi)$,
then since
\[
\left\{\frac{e^{i\ell s}}{\sqrt{2\pi}}\Big|_{[0,2\pi)}\colon \ \ell\in
\mathbb{Z}\right\}
\]
is an orthonormal basis for $L^2[0,2\pi]$ and the exponentials $e^{i\ell s}$ are
$2\pi$-invariant, as in the case of Shannon's wavelet it follows that
\[
\left\{\frac{e^{i\ell s}}{\sqrt{2\pi}}\Big|_E\colon \ \ell\in \mathbb{Z}\right\}
\]
is an orthonormal basis for $L^2(E)$. Also, if $E$ is $2\pi$-translation
congruent to $[0,2\pi)$, then since
\[
\{[0,2\pi) + 2\pi n\colon \ n\in \mathbb{Z}\}
\]
is a measurable partition of $\mathbb{R}$, so is
\[
\{E + 2\pi n\colon \ n\in \mathbb{Z}\}.
\]
These arguments can be reversed.

We say that a measurable subset $G\subset \mathbb{R}$ is a 2-{\em dilation
generator\/} of a {\em partition\/} of $\mathbb{R}$ if the sets
\begin{equation}\label{eq41}
2^nG := \{2^ns\colon \ s\in G\},\qquad n\in \mathbb{Z}
\end{equation}
are disjoint and $\mathbb{R}\backslash \cup_n 2^nG$ is a null set. Also, we say
that $E\subset \mathbb{R}$ is a $2\pi$-{\em translation generator of a
partition\/} of $\mathbb{R}$ if the sets
\begin{equation}\label{eq42}
E + 2n\pi := \{s + 2 n\pi\colon \ s\in E\},\qquad n\in \mathbb{Z},
\end{equation}
are disjoint and $\mathbb{R}\backslash \cup_n (E+2n\pi)$ is a null set.

\begin{lemma}\label{lem8}
A measurable set $E\subseteq \mathbb{R}$ is a $2\pi$-translation generator of a
partition of $\mathbb{R}$ if and only if, modulo a null set, $E$ is translation
congruent to $[0,2\pi)$ modulo $2\pi$. Also, a measurable set $G\subseteq
\mathbb{R}$ is a 2-dilation generator of a partition of $\mathbb{R}$ if and only
if, modulo a null set, $G$ is a dilation congruent modulo 2 to the set $[-2\pi,
-\pi) \cup [\pi,2\pi)$.
\end{lemma}

\subsubsection{A Criterion}
The following is a useful criterion for wavelet sets. It was
published independently by Dai--Larson in [11] and by Fang--Wang in
[17] at about the same time in December, 1994. In fact, it is
amusing that the two papers had been submitted within two days of
each other; only much later did we even learn of each others work
and of this incredible timing.

\begin{proposition}\label{pro9}
Let $E\subseteq\mathbb{R}$ be a measurable set. Then $E$ is a wavelet set if and
only if $E$ is both a 2-dilation generator of a partition (modulo null sets) of
$\mathbb{R}$ and a $2\pi$-translation generator of a partition (modulo null
sets) of $\mathbb{R}$. Equivalently, $E$ is a wavelet set if and only if $E$ is
both translation congruent to $[0,2\pi)$ modulo $2\pi$ and dilation congruent to
$[-2\pi,-\pi) \cup [\pi,2\pi)$ modulo 2.
\end{proposition}

Note that a set is $2\pi$-translation congruent to $[0,2\pi)$ iff it is
$2\pi$-translation congruent to $[-2\pi, \pi)\cup [\pi,2\pi)$. So the last
sentence of Proposition \ref{pro9} can be stated:\ A measurable set $E$ is a
wavelet set if and only if it is both $2\pi$-translation and 2-dilation
congruent to the Littlewood--Paley set $[-2\pi, -\pi)\cup [\pi,2
\pi)$.

\subsection{Phases}
If $E$ is a wavelet set, and if $f(s)$ is any function with support
$E$ which has constant modulus $\frac1{\sqrt{2\pi}}$ on $E$, then
$\mathcal{F}^{-1}(f)$ is a wavelet. Indeed, by Lemma \ref{lem7}
$\{\widehat T^\ell f\colon \ \in \mathbb{Z}\}$ is an orthonormal
basis for $L^2(E)$, and since the sets $2^nE$ partition
$\mathbb{R}$, so $L^2(E)$ is a complete wandering subspace for
$\widehat D$, it follows that $\{\widehat D^n\widehat T^\ell
f\colon\ n,\ell\in \mathbb{Z}\}$ must be an orthonormal basis for
$L^2(\mathbb{R})$, as required. In [17, 37, 38] the term MSF-wavelet
includes this type of wavelet. So MSF-wavelets can have arbitrary
phase and $s$-elementary wavelets have phase 0. {\em Every\/} phase
is attainable in the sense of chapter 3 for an MSF or $s$-elementary
wavelet.

\subsubsection{Some Examples of One-Dimensional Wavelet Sets}
It is usually easy to determine, using the dilation-translation criteria, in
Proposition \ref{pro9}, whether a given finite union of intervals is a wavelet
set. In fact, to verify that a given ``candidate'' set $E$ is a wavelet set, it
is clear from the above discussion and criteria that it suffices to do two
things.

\begin{quote}
(1)~~Show, by appropriate partitioning, that $E$ is 2-dilation-congruent to a
set of the form $[-2a,-a)\cup [b,2b)$ for some $a,b>0$.
\end{quote}
and
\begin{quote}
(2)~~Show, by appropriate partitioning, that $E$ is $2\pi$-translation-congruent
to a set of the form $[c,c+2\pi)$ for some real number $c$.
\end{quote}

On the other hand, wavelet sets suitable for testing hypotheses, can be quite
difficult to construct. There are very few ``recipes'' for wavelet sets, as it
were. Many families of such sets have been constructed for reasons including
perspective, experimentation, testing hypotheses, etc., including perhaps the
pure enjoyment of doing the computations -- which are somewhat ``puzzle-like''
in nature. In working with the theory it is nice (and in fact necessary) to have
a large supply of wavelets on hand that permit relatively simple analysis.

For this reason we take the opportunity here to present for the
reader a collection of such sets, mainly taken from [11], leaving
most of the ``fun'' in verifying that they are indeed wavelet sets
to the reader.

We refer the reader to [12] for a proof of the existence of wavelet
sets in $\mathbb{R}^{(n)}$, and a proof that there are sufficiently
many to generate the Borel structure of $\mathbb{R}^{(n)}$. These
results are true for arbitrary expansive dilation factors. Some
concrete examples in the plane were subsequently obtained by Soardi
and Weiland, and others were obtained by Gu and Speegle. Two had
also been obtained by Dai for inclusion in the revised concluding
remarks section of our Memoir [11].

In these examples we will usually write intervals as half-open intervals
$[\cdot,~)$ because it is easier to verify the translation and dilation
congruency relations (1) and (2) above when wavelet sets are written thus, even
though in actuality the relations need only hold modulo null sets.

(i)~~As mentioned above, an example due to Journe of a wavelet which admits no
multiresolution analysis is the $s$-elementary wavelet with wavelet set
\[
\left[-\frac{32\pi}7, -4\pi\right)\cup \left[-\pi, \frac{4\pi}7\right)\cup
\left[\frac{4\pi}7, \pi\right) \cup \left[4\pi, \frac{32\pi}7\right).
\]
To see that this satisfies the criteria, label these intervals, in order, as
$J_1, J_2, J_3, J_4$ and write $J =\cup J_i$. Then
\[
J_1\cup 4J_2 \cup 4J_3\cup J_4 = \left[-\frac{32\pi}7, -\frac{16\pi}7\right)
\cup \left[\frac{16\pi}7, \frac{32\pi}7\right).\]
This has the form $[-2a,a)\cup [b,2b)$ so is a 2-dilation generator of a
partition of $\mathbb{R}\backslash\{0\}$. Then also observe that
\[
\{J_1 + 6\pi, J_2 +2\pi, J_3, J_4-4\pi\}
\]
is a partition of $[0,2\pi)$.

(ii)~~The Shannon (or Littlewood--Paley) set can be generalized. For
any $-\pi < \alpha < \pi$, the set
\[
E_\alpha = [-2\pi + 2\alpha, -\pi + \alpha) \cup [\pi + \alpha, 2\pi +2\alpha)
\]
is a wavelet set. Indeed, it is clearly a 2-dilation generator of a partition of
$\mathbb{R}\backslash\{0\}$, and to see that it satisfies the translation
congruency criterion for $-\pi < \alpha \le 0$ (the case $0<\alpha<\pi$ is
analogous) just observe that
\[
\{[-2\pi + 2\alpha, 2\pi) + 4\pi, [-2\pi, -\pi+\alpha) + 2\pi, [\pi +
\alpha,2\pi + 2\alpha)\}
\]
is a partition of $[0,2\pi)$. It is clear that $\psi_{E_\alpha}$ is then a
continuous (in $L^2(\mathbb{R})$-norm) path of $s$-elementary wavelets. Note
that
\[
\lim_{\alpha\to\pi} \widehat\psi_{E_\alpha} = \frac1{\sqrt{2\pi}} \chi_{[2\pi,
4\pi)}.
\]
This is {\em not\/} the Fourier transform of a wavelet because the set $[2\pi,
4\pi)$ is not a 2-dilation generator of a partition of
$\mathbb{R}\backslash\{0\}$. So
\[
\lim_{\alpha\to\pi} \psi_{E_\alpha}
\]
is not an orthogonal wavelet. (It is what is known as a Hardy wavelet because it
generates an orthonormal basis for $H^2(\mathbb{R})$ under dilation and
translation.) This example demonstrates that $\mathcal{W}(D,T)$ is {\em not\/}
closed in $L^2(\mathbb{R})$.

(iii)~~Journe's example above can be extended to a path. For $-\frac\pi7 \le
\beta\le \frac\pi7$ the set
\[
J_\beta = \left[-\frac{32\pi}7, -4\pi + 4\beta\right) \cup \left[-\pi +\beta,
-\frac{4\pi}7\right) \cup \left[\frac{4\pi}7, \pi+\beta\right)\cup \left[4 \pi +
4\beta, 4\pi + \frac{4\pi}7\right)
\]
is a wavelet set. The same argument in (i) establishes dilation congruency. For
translation, the argument in (i) shows congruency to $[4\beta, 2\pi+4\beta)$
which is in turn congruent to $[0,2\pi)$ as required. Observe that here, as
opposed to in (ii) above, the limit of $\psi_{J_\beta}$ as $\beta$ approaches
the boundary point $\frac\pi7$ {\em is\/} a wavelet. Its wavelet set is a union
of 3 disjoint intervals.

(iv)~~Let $A\subseteq [\pi, \frac{3\pi}2)$ be an arbitrary measurable subset.
Then there is a wavelet set $W$, such that $W\cap [\pi, \frac{3\pi}2)=A$. For
the
construction, let
\begin{align*}
B &= [2\pi, 3\pi)\backslash 2A,\\
C &= \left[-\pi, -\frac\pi2\right)\backslash (A-2\pi)\\
\text{and}\quad D &= 2A-4\pi.
\end{align*}
Let
\[
W = \left[\frac{3\pi}2, 2\pi\right)\cup A\cup B\cup C\cup D.
\]
We have $W\cap [\pi, \frac{3\pi}2) =A$. Observe that the sets
$[\frac{3\pi}2,2\pi)$, $A,B,C,D$, are disjoint. Also observe that the sets
\[
\left[\frac{3\pi}2,2\pi\right), A, \frac12 B, 2C, D,
\]
are disjoint and have union $[-2\pi,-\pi)\cup [\pi,2\pi)$. In addition, observe
that the sets
\[
\left[\frac{3\pi}2,2\pi\right), A,B-2\pi, C + 2\pi, D+2\pi,\]
are disjoint and have union $[0,2\pi)$. Hence $W$ is a wavelet set.

(v)~~Wavelet sets for arbitrary (not necessarily integral) dilation factors
other then 2 exist. For instance, if $d\ge 2$ is arbitrary, let
\begin{align*}
A &= \left[-\frac{2d\pi}{d+1}, -\frac{2\pi}{d+1}\right),\\
B &= \left[\frac{2\pi}{d^2-1}, \frac{2\pi}{d+1}\right),\\
C &= \left[\frac{2d\pi}{d+1}, \frac{2d^2\pi}{d^2-1}\right)
\end{align*}
and let $G = A \cup B\cup C$. Then $G$ is $d$-wavelet set. To see this, note
that
$\{A+2\pi,B,C\}$ is a partition of an interval of length $2\pi$. So $G$ is
$2\pi$-translation-congruent to $[0,2\pi)$. Also, $\{A,B,d^{-1}C\}$ is a
partition of the set $[-d\alpha, -\alpha)\cup [\beta,d\beta)$ for $\alpha =
\frac{2\pi}{d^2-1}$, and $\beta = \frac{2\pi}{d^2-1}$, so from this form it
follows that $\{d^nG\colon \ n\in \mathbb{Z}\}$ is a partition of
$\mathbb{R}\backslash\{0\}$. Hence if $\psi :=
\mathcal{F}^{-1}(\frac1{\sqrt{2\pi}} \chi_G)$, it follows that
$\{d^{\frac{n}2}\psi(d^nt-\ell)\colon \ n,\ell\in \mathbb{Z}\}$ is orthonormal
basis for $L^2(\mathbb{R})$, as required.

\subsection{Operator-Theoretic Interpolation of Wavelets: The Special Case of Wavelet Sets}

Let $E, F$ be a pair of wavelet sets.  Then for (a.e.) $x \in E$
there is a unique $y \in F$ such that $x - y \in 2\pi\mathbb{Z}$.
This is the $\emph{translation congruence}$ property of wavelet
sets.  Also, for (a.e.) $x \in E$ there is a unique $z \in F$ such
that $\frac{x}{z}$ is an integral power of $2$. This is the
\emph{dilation congruence} property of wavelet sets. (See section
2.5.6.)

There is a natural \emph{closed-form algorithm} for the
$\emph{interpolation unitary}$ $V_{\psi_E}^{\psi_F}$ which maps the
wavelet basis for  $\widehat{\psi}_E $ to the wavelet basis for
$\widehat{\psi}_F$.  Indeed, using both the translation and dilation
congruence properties of $\{E,F\}$, one can explicitly compute a
(unique) measure-preserving transformation $\sigma := {\sigma}_E^F$
mapping $\mathbb{R}$ onto $\mathbb{R}$ which has the property that
$V_{\psi_E}^{\psi_F}$ is identical with the \emph{composition
operator} defined by:
$$ f \mapsto f \circ {\sigma}^{-1} $$
for all $f \in L^2(\mathbb{R})$. With this formulation, compositions
of the maps $\sigma$ between different pairs of wavelet sets are not
difficult to compute, and thus products of the corresponding
interpolation unitaries can be computed in terms of them.

\subsubsection{The Interpolation Map $\sigma$}
Let $E$ and $F$ be arbitrary wavelet sets. Let $\sigma\colon \ E\to
F$ be the 1-1, onto map implementing the $2\pi$-translation
congruence. Since $E$ and $F$ both generated partitions of
$\mathbb{R}\backslash\{0\}$ under dilation by powers of 2, we may
extend $\sigma$ to a 1-1 map of $\mathbb{R}$ onto $\mathbb{R}$ by
defining $\sigma(0)=0$, and
\begin{equation}\label{eq44}
\sigma(s) = 2^n\sigma(2^{-n}s) \quad \text{for}\quad s\in 2^nE,\quad
n\in \mathbb{Z}.
\end{equation}
We adopt the notation $\sigma^F_E$ for this, and call it the {\em
interpolation map\/} for the ordered pair $(E,F)$.

\begin{lemma}\label{lem5.1}
In the above notation, $\sigma^F_E$ is a measure-preserving
transformation from $\mathbb{R}$ onto $\mathbb{R}$.
\end{lemma}

\begin{proof}
Let $\sigma := \sigma^F_E$. Let $\Omega\subseteq \mathbb{R}$ be a
measurable set. Let $\Omega_n= \Omega \cap 2^nE$, $n\ni\mathbb{Z}$,
and let $E_n = 2^{-n} \Omega_n\subseteq E$. Then $\{\Omega_n\}$ is a
partition of $\Omega$, and we have $m(\sigma(E_n)) = m(E_n)$ because
the restriction of $\sigma$ to $E$ is measure-preserving. So
\begin{align*}
m(\sigma(\Omega)) &= \sum_n m(\sigma(\Omega_n)) = \sum_n m(2^n \sigma(E_n))\\
&= \sum_n 2^nm(\sigma(E_n)) = \sum_n 2^nm(E_n)\\
&= \sum_n m(2^nE_n) =\sum_n m(\Omega_n) = m(\Omega).
\end{align*}
\end{proof}

A function $f\colon \mathbb{R}\to \mathbb{R}$ is called 2-{\em
homogeneous\/} if $f(2s) = 2f(s)$ for all $s\in \mathbb{R}$.
Equivalently, $f$ is 2-homogeneous iff $f(2^ns) = 2^nf(s)$,
$s\in\mathbb{R}$, $n\in\mathbb{Z}$. Such a function is completely
determined by its values on any subset of $\mathbb{R}$ which
generates a partition of $\mathbb{R}\backslash\{0\}$ by 2-dilation.
So $\sigma^F_E$ is the (unique) 2-homogeneous extension of the
$2\pi$-transition congruence $E\to F$. The set of all 2-homogeneous
measure-preserving transformations of $\mathbb{R}$ clearly forms a
group under composition. Also, the composition of a
2-dilation-periodic function $f$ with a 2-homogeneous function $g$
is (in either order) 2-dilation periodic. We have $f(g(2s)) =
f(2g(s)) = f(g(s))$ and  $g(f(2s)) = g(f(s))$. These facts will be
useful.

\subsubsection{An Algorithm For The Interpolation Unitary}

 Now let
\begin{equation}\label{eq45}
U^F_E := U_{\sigma^F_E},
\end{equation}
where if $\sigma$ is any measure-preserving transformation of
$\mathbb{R}$ then $U_\sigma$ denotes the composition operator
defined by $U_\sigma f = f\circ\sigma^{-1}$, $f\in L^2(\mathbb{R})$.
Clearly $(\sigma^F_E)^{-1} = \sigma^E_F$ and $(U^F_E)^* = U^E_F$. We
have $U^F_E\widehat\psi_E = \widehat\psi_F$ since $\sigma^F_E(E)=F$.
That is,
\[
U^F_E\widehat\psi_E = \widehat\psi_E \circ\sigma^E_F =
\frac1{\sqrt{2\pi}} \chi_{_E} \circ \sigma^E_F = \frac1{\sqrt{2\pi}}
\chi_{_F} = \widehat\psi_F.
\]

\begin{proposition}\label{pro5.2}
Let $E$ and $F$ be arbitrary wavelet sets. Then $U^F_E\in
\mathcal{C}_{\widehat\psi_E}(\widehat D,\widehat T)$.  Hence
$\mathcal{F}^{-1} U^F_E \mathcal{F}$ is the interpolation unitary
for the ordered pair $(\psi_E,\psi_F)$.
\end{proposition}

\begin{proof}
Write $\sigma = \sigma^F_E$ and $U_\sigma = U^F_E$. We have
$U_\sigma\widehat\psi_E = \widehat\psi_F$ since $\sigma(E) = F$. We
must show
\[
U_\sigma\widehat D^n \widehat T^l \widehat\psi_E = \widehat
D^n\widehat T^l U_\sigma \widehat\psi_E,\quad n,l\in \mathbb{Z}.
\]
We have
\begin{align*}
(U_\sigma\widehat D^n\widehat T^l \widehat\psi_E)(s) &=
(U_\sigma\widehat D^n
e^{-ils} \widehat\psi_E)(s)\\
&= U_\sigma2^{-\frac{n}2} e^{-il2^{-n}s} \widehat\psi_E(2^{-n}s)\\
&= 2^{-\frac{n}2} e^{-il2^{-n}\sigma^{-1}(s)}
\widehat\psi_E(2^{-n}\sigma^{-1}(s))\\
&= 2^{-\frac{n}2} e^{-il\sigma^{-1}(2^{-n}s)}
\widehat\psi_E(\sigma^{-1}(2^{-n}s))\\
&= 2^{-\frac{n}2} e^{-il\sigma^{-1}(2^{-n}s)} \widehat\psi(2^{-n}s).
\end{align*}
This last term is nonzero iff $2^{-n}s\in F$, in which case
$\sigma^{-1}(2^{-n}s) = \sigma^E_F(2^{-n}s)$ $=2^{-n}s+2\pi k$ for
some $k\in \mathbb{Z}$ since $\sigma^E_F$ is a
$2\pi$-translation-congruence on $F$. It follows that
$e^{-il\sigma^{-2}(2^{-n}s)} = e^{-il2^{-n}s}$. Hence we have
\begin{align*}
(U_\sigma\widehat D^n\widehat T^l\widehat\psi_E)(s) &=
2^{-\frac{n}2}
e^{-ils^{-2n}s} \widehat \psi_F(2^{-n}s)\\
&= (\widehat D^n\widehat T^l\widehat\psi_F)(s)\\
&= (\widehat D^n\widehat T^lU_\sigma\widehat \psi_E)(s).
\end{align*}
We have shown $U^F_E \in \mathcal{C}_{\widehat\psi_E}(\widehat
D,\widehat T)$. Since $U^F_E\widehat\psi_E = \widehat\psi_F$, the
uniqueness part of Proposition ~1 shows that
$\mathcal{F}^{-1}U^F_E\mathcal{F}$ must be the interpolation unitary
for $(\psi_E,\psi_F)$.
\end{proof}

\subsection{The Interpolation Unitary Normalizes The Commutant}

\begin{proposition}\label{pro5.3}
Let $E$ and $F$ be arbitrary wavelet sets. Then the interpolation
unitary for the ordered pair $(\psi_E,\psi_F)$ normalizes
$\{D,T\}'$.
\end{proposition}

\begin{proof}
By Proposition \ref{pro5.2} we may work with $U^F_E$ in the Fourier
transform domain. By Theorem 6, the generic element of $\{\widehat
D,\widehat T\}'$ has the form $M_h$ for some 2-dilation-periodic
function $h\in L^\infty(\mathbb{R})$. Write $\sigma = \sigma^F_E$
and $U_\sigma= U^F_E$. Then
\begin{equation}\label{eq46}
U^{-1}_\sigma M_hU_\sigma= M_{h\circ\sigma^{-1}}.
\end{equation}
So since the composition of a 2-dilation-periodic function with a
2-homogeneous function is 2-dilation-periodic, the proof is
complete.
\end{proof}

\subsubsection{$\mathcal{C}_\psi(D,T)$ is Nonabelian}

It can also be shown ([11, Theorem 5.2 (iii)]) that if $E,F$ are
wavelet sets with $E\ne F$ then $U^F_E$ is not contained in the
double commutant $\{\widehat D,\widehat T\}''$. So since $U^F_E$ and
$\{\widehat D,\widehat T\}'$ are both contained in the local
commutant of $\mathcal{U}_{\widehat D,\widehat T}$ at
$\widehat\psi_E$, this proves that
$\mathcal{C}_{\widehat\psi_E}(\widehat D,\widehat T)$ is nonabelian.
In fact (see [11, Proposition 1.8]) this can be used to show that
$\mathcal{C}_\psi(D,T)$ is nonabelian for every wavelet $\psi$. We
suspected this, but we could not prove it until we discovered the
``right'' way of doing the needed computation using $s$-elementary
wavelets.

The above shows that a pair $(E,F)$ of wavelets sets (or, rather,
their corresponding $s$-elementary wavelets) admits
operator-theoretic interpolation if and only if Group$\{U^F_E\}$ is
contained in the local commutant
$\mathcal{C}_{\widehat\psi_E}(\widehat D,\widehat T)$, since the
requirement that $U^F_E$ normalizes $\{\widehat D,\widehat T\}'$ is
automatically satisfied. It is easy to see that this is equivalent
to the condition that for each $n\in\mathbb{Z}$, $\sigma^n$ is a
$2\pi$-congruence of $E$ in the sense that
$(\sigma^n(s)-s)/2\pi\in\mathbb{Z}$ for all $s\in E$, which in turn
implies that $\sigma^n(E)$ is a wavelet set for all $n$. Here
$\sigma = \sigma^F_E$. This property hold trivially if $\sigma$ is
{\em involutive\/} (i.e.\ $\sigma^2=$ identity).

\subsubsection{The Coefficient Criterion}
In cases where ``torsion'' is present, so $(\sigma^F_E)^k$ is the
identity map for some finite integer $k$, the von Neumann algebra
generated by $\{\widehat D,\widehat T\}'$ and $U :=U^F_E$ has the
simple form
\[
\left\{\sum^k_{n=0} M_{h_n}U^n\colon \ h_n\in L^\infty(\mathbb{R})
\text{ with } h_n(2s) = h_n(s),\quad s\in \mathbb{R}\right\},
\]
and so each member of this ``interpolated'' family of wavelets has
the form
\begin{equation}\label{eq47}
\frac1{\sqrt{2\pi}} \sum^k_{n=0} h_n(s) \chi_{\sigma^n(E)}
\end{equation}
for 2-dilation periodic ``coefficient'' functions $\{h_n(s)\}$ which
satisfy the necessary and sufficient condition that the operator
\begin{equation}\label{eq48}
\sum^k_{n=0} M_{h_n }U^n
\end{equation}
is unitary.

A standard computation shows that the map $\theta$ sending $\sum^k_0
M_{h_n}U^n$ to the $k\times k$ function matrix $(h_{ij})$ given by
\begin{equation}\label{eq49}
h_{ij} = h_{\alpha(i,j)}\circ \sigma^{-i+1}
\end{equation}
where $\alpha(i,j) = (i+1)$ modulo $k$, is a $*$-isomorphism. This
matricial algebra is the cross-product of $\{D,T\}'$ by the
$*$-automorphism $ad(U^F_E)$ corresponding to conjugation with
$U^F_E$. For instance,  if $k=3$ then $\theta$ maps
\[
M_{h_1} + M_{h_2} U^F_E + M_{h_3}(U^F_E)^2
\]
to
\begin{equation}\label{eq50}
\left(\begin{matrix}
h_1&h_2&h_3\\
h_3\circ\sigma^{-1}&h_1\circ\sigma^{-1}&h_2\circ\sigma^{-1}\\
h_2\circ\sigma^{-2}&h_3\circ\sigma^{-2}&h_1\circ\sigma^{-2}
\end{matrix}\right).
\end{equation}
This shows that $\sum^k_0 M_{h_n}U^n$ is a unitary operator iff the
scalar matrix $(h_{ij})(s)$ is unitary for almost all $s\in
\mathbb{R}$. Unitarity of this matrix-valued function is called the
{\em Coefficient Criterion\/} in [11], and the functions $h_i$ are
called the interpolation coefficients. This leads to formulas for
families of wavelets which are new to wavelet theory.

\subsection{Interpolation Pairs of Wavelet Sets}

 For many
interesting cases of note, the interpolation map  $\sigma^F_E$ will
in fact be an \emph{involution} of $\mathbb{R}$ (i.e. $\sigma \circ
\sigma = id$, where $\sigma := \sigma^F_E$, and where $id$ denotes
the identity map). So torsion \emph{will} be present, as in the
above section, and it will be present in an essentially simple form.
The corresponding interpolation unitary will be a \emph{symmetry} in
this case (i.e. a selfadjoint unitary operator with square $I$).

It is curious to note that verifying a simple operator equation $U^2
= I$ directly by matricial computation can be extremely difficult.
It is much more computationally feasible to verify an equation such
as this by pointwise (a.e.) verifying explicitly the relation
$\sigma \circ \sigma = id$ for the interpolation map.  In [11] we
gave a number of examples of interpolation pairs of wavelet sets. We
give below a collection of examples that has not been previously
published: Every pair sets from the Journe family is an
interpolation pair.

\subsection{Journe Family Interpolation Pairs}

Consider the parameterized path of \emph{generalized Journe} wavelet
sets given in [11, Example 4.5(iii)].  We have
$$ J_{\beta} = \left[-\frac{32\pi}{7} , -4\pi - 4\beta\right) \cup \left[-\pi +
\beta
, -\frac{4\pi}{7}\right) \cup\left[\frac{4\pi}{7}, \pi + \beta\right) \cup
\left[4\pi +
4\beta, 4\pi + \frac{4\pi}{7}\right)$$
where the set of parameters $\beta$ ranges $-\frac{\pi}{7} \leq \beta
\leq \frac{\pi}{7}$.

\begin{proposition} Every pair $(J_{\beta_1}, J_{\beta_2})$ is an
interpolation pair.
\end{proposition}
\begin{proof}
Let $\beta_1, \beta_2 \in \left[-\frac{\pi}{7}, \frac{\pi}{7}\right)$ with
$\beta_1
< \beta_2$. Write $\sigma = \sigma_{J_{\beta_2}}^{J_{\beta_1}} .$ We
need to show that
\begin{equation}{\sigma}^2(x)=x \tag{*}
\end{equation} for all $x \in \mathbb{R}$. Since $\sigma$ is
2-homogeneous, it suffices to verify (*) only for $x \in
J_{\beta_1}$. For $x \in J_{\beta_1} \cap J_{\beta_2}$ we have
$\sigma (x) = x$, hence ${\sigma}^2(x)=x$.  So we only need to check
(*) for $x \in (J_{\beta_1} \backslash J_{\beta_2})$. We have
$$J_{\beta_1} \backslash J_{\beta_2} = [-\pi + \beta_1 , -\pi +
\beta_2) \cup [4\pi + 4\beta_1 , 4\pi + 4\beta_2) .$$ It is useful
to also write
$$ J_{\beta_2} \backslash J_{\beta_1} = [-4\pi + 4\beta_1 , -4\pi +
4\beta_2) \cup [\pi + \beta_1 , \pi + \beta_2).$$

On $[-\pi + \beta_1 , -\pi + \beta_2 )$ we have $\sigma (x) = x +
2\pi$, which lies in $[\pi + \beta_1 , \pi + \beta_2)$.  If we
multiply this by $4$, we obtain
$4\sigma(x)\in[4\pi+4\beta_1,4\pi+4\beta_2) \subset J_{\beta_1}$.
And on $[4\pi + 4\beta_1, 4\pi + 4\beta_2)$ we clearly have $\sigma
(x) = x - 8\pi$, which lies in $[-4\pi + 4\beta_1 , -4\pi +
4\beta_2)$.

So for $x \in [-\pi + \beta_1 , -\pi + \beta_2)$ we have
$${\sigma}^2(x) = \sigma (\sigma (x) ) = \frac14 \sigma (4\sigma
(x)) = \frac14 [4\sigma(x) - 8\pi] = \sigma(x) -2\pi = x + 2\pi
- 2\pi = x.$$

On $[4\pi + 4\beta_1, 4\pi +4\beta_2)$ we have $\sigma(x) = x -
8\pi$, which lies in $[-4\pi + 4\beta_1 , -4\pi + 4\beta_2)$. So $\frac14 \sigma
(x) \in [-\pi + \beta_1, -\pi+\beta_2)$. Hence
$$\sigma\left(\frac14\sigma(x)\right) = \frac14\sigma(x) + 2\pi$$
and
thus
$${\sigma}^2(x) = 4\sigma\left(\frac14\sigma(x)\right) =
4\left[\frac14\sigma(x) + 2\pi\right] = \sigma(x) + 8\pi = x - 8\pi + 8\pi =
x$$ as required.

We have shown that for all  $x \in J_{\beta_1}$ we have
${\sigma}^2(x) = x$.  This proves that $(J_{\beta_1}, J_{\beta_2})$
is an interpolation pair.

\end{proof}

\section{Unitary Systems and Frames}

In [33] we developed an operator-theoretic approach to discrete
frame theory (i.e.\ frame sequences, as opposed to continuous frame
transforms) on a separable Hilbert space. We then applied it to an
investigation of frame vectors for unitary systems, frame wavelets
and group representations. The starting-point idea, which is pretty
simple-minded in fact, is to realize any frame sequence for a
Hilbert space $H$ as a compression of a Riesz basis for a larger
Hilbert space. In other words, a frame is a sequence of vectors in a
Hilbert space which {\em dilates}, (in the operator-theoretic or
geometric sense, as opposed to the function-theoretic sense of
multiplication of the independent variable of a function by a
dilation constant), or {\em extends}, to a (Riesz) basis for a
larger space. From this idea much can be developed, and some new
perspective can be given to certain concepts that have been used in
engineering circles for many years. See section \ref{SA} below.

\subsection{Basics on Frames}
Let $H$ be a separable complex Hilbert space. Let $B(H)$ denote the
algebra of all bounded linear operators on $H$. Let $\mathbb{N}$
denote the natural numbers, and $\mathbb{Z}$ the integers. We will
use $\mathbb{J}$ to denote a generic countable (or finite) index set
such as $\mathbb{Z}, \mathbb{N}, \mathbb{Z}^{(2)}$,
$\mathbb{N}\cup\mathbb{N}$ etc.

A sequence $\{x_j\colon \ j\in \mathbb{N}\}$ of vectors in $H$ is
called a {\em frame\/} if there are constants $A,B>0$ such that
\[
A\|x\|^2 \le \sum_j |\langle x,x_j\rangle|^2 \le B\|x\|^2
\]
for all $x\in H$. The optimal constant (maximal for $A$ and minimal
for $B$) are called the {\emph frame bounds}. The frame $\{x_j\}$ is
called a {\em tight frame\/} if $A=B$, and is called {\emph
Parseval\/} if $A=B=1$. (Originally, in [33] and a in number of
subsequent papers, the term \emph{normalized tight frame} was used
for this.  However, this term had also been applied by Benedetto and
Ficus [5] for another concept: a tight frame of unit vectors; what
we now call a uniform tight frame, or spherical frame. So, after all
parties involved, the name \emph{Parseval} was adopted. It makes a
lot of sense, because a Parseval frame is precisely a frame which
satisfies Parseval's identity.) A sequence $\{x_j\}$ is defined to
be a {\emph Riesz basis\/} if it is a frame and is also a basis for
$H$ in the sense that for each $x\in H$ there is a {\em unique\/}
sequence $\{\alpha_j\}$ in $\mathbb{C}$ such that $x = \sum
\alpha_jx_j$ with the convergence being in norm. We note that a
Riesz basis is also defined to be  basis which is obtained from an
orthonormal basis by applying a bounded linear invertible operator.
This is equivalent to the first definition. It should be noted that
in Hilbert spaces the Riesz bases are precisely the bounded
unconditional bases. We will say that frames $\{x_j\colon \ j\in
\mathbb{J}\}$ and $\{y_j\colon \ j\in \mathbb{J}\}$ on Hilbert
spaces $H,K$, respectively, are {\em unitarily equivalent\/} if
there is a unitary $U\colon \ H\to K$ such that $Ux_j = y_j$ for all
$j\in \mathbb{J}$. We will say that they are {\em similar\/} (or
{\em isomorphic\/}) if there is a bounded linear invertible operator
$T\colon \ H\to K$ such that $Tx_j=y_j$ for all $j\in \mathbb{J}$.

\begin{example}\label{exmA}
Let $K=L^2(\mathbb{T})$ where $\mathbb{T}$ is the unit circle and
measure is normalized Lebesgue measure, and let $\{e^{ins}\colon \
n\in \mathbb{Z}\}$ be the standard orthonormal basis for
$L^2(\mathbb{T})$. If $E\subseteq \mathbb{T}$ is any measurable
subset then $\{e^{ins}|_E\colon \ n\in \mathbb{Z}\}$ is a Parseval
frame for $L^2(E)$. This can be viewed as obtained from the single
vector $\chi_{_E}$ by applying all integral powers of the (unitary)
multiplication operator $M_{e^{is}}$. It turns out that these are
all (for different $E$) unitarily {\em inequivalent}. This is an
example of a Parseval frame which is generated by the action of a
unitary group on a single vector. This can be compared with the
definition of a \em{frame wavelet}. (As one might expect, a single
function $\psi$ in $L^2(\mathbb{R})$ which generates a {\em frame\/}
for $L^2(\mathbb{R})$ under the action of $\mathcal{U}_{D,T}$ is
called a {\em frame-wavelet}.)
\end{example}

\subsection{Dilation of Frames: The Discrete Version of Naimark's
Theorem}\label{SA}

Now let $\{x_n\}_{n\in \mathbb{J}}$ be a Parseval frame and let
$\theta\colon \ H\to K := l^2(\mathbb{J})$ be the usual {\em
analysis operator\/} (this was called the {\em frame transform in
[HL]\/}) defined by $\theta(x) := (\langle x,x_n\rangle)_{n\in
\mathbb{J}}$. This is obviously an isometry. Let $P$ be the
orthogonal projection from $K$ onto $\theta(H)$. Denote the standard
orthonormal basis for $l^2(\mathbb{J})$ by $\{e_j\colon \ j\in
\mathbb{J}\}$. For any $m\in \mathbb{J}$, we have
\begin{align*}
\langle\theta(x_m),Pe_n\rangle &= \langle P\theta(x_m),e_n\rangle =
\langle\theta (x_m),e_n\rangle\\
&= \langle x_m,x_n\rangle = \langle\theta(x_m), \theta(x_m)\rangle.
\end{align*}
It follows easily that $\theta(x_n) = Pe_n$, $n\in\mathbb{J}$.
Identifying $H$ with $\theta(H)$, this shows indeed that every
Parseval frame can be realized by compressing an orthonormal basis,
as claimed earlier.

 This can actually be viewed as a special case
(probably the simplest possible special case) of an old theorem of
Naimark concerning operator algebras and dilation of positive
operator valued measures to projection valued measures.  The
connection between Naimark's theorem and the dilation result for
Parseval frames, and that the latter can be viewed as a special case
of the former, was pointed out to me by Chandler Davis and Dick
Kadison in a conference (COSY-1999: The Canadian Operator Algebra
Symposium, Prince Edward Island, May 1999).

\subsection{Complements of Frames}
It is useful to note that $P$ will equal $I$ iff $\{x_n\}$ is a basis. Indeed,
if $P\ne I$, then choose $z\ne 0$, $z\in (I-P)K$, and write $z=\sum \alpha_ne_n$
for some sequence $\alpha_n\in\mathbb{C}$. Then $0 = Pz = \sum
\alpha_n\theta(x_n)$, and not all the scalars $\alpha_n$ are zero. Hence
$\{x_n\}$ is not topologically linearly independent so cannot even be a Schauder
basis. On the other hand if $P=I$ then $\{x_n\}$ is obviously an orthonormal
basis.

Suppose $\{x_n\}_{n\in\mathbb{J}}$ is a Parseval frame for $H$, and
let $\theta,P,K,e_n$ be as above. Let $M = (I-P)K$. Then $y_n :=
(I-P)e_n$ is a Parseval frame on $M$ which is {\em complementary\/}
to $\{x_n\}$ in the sense that the inner direct sum $\{x_n\oplus
y_n\colon \ n\in\mathbb{J}\}$ is an orthonormal basis for the direct
sum Hilbert space $H\oplus M$. Moreover there is uniqueness:\ The
extension of a tight frame to an orthonormal basis described in the
above paragraph is unique up to unitary equivalence. That is if $N$
is another Hilbert space and $\{z_n\}$ is a tight frame for $N$ such
that $\{x_n\oplus z_n\colon \ n\in \mathbb{J}\}$ is an orthonormal
basis for $H\oplus N$, then there is a unitary transformation $U$
mapping $M$ onto $N$ such that $Uy_n = z_n$ for all $n$. In
particular, $\dim M = \dim N$.

If $\{x_j\}$ is a Parseval frame, we will call any Parseval frame
$\{z_j\}$ such that $\{x_j\oplus z_j\}$ is an orthonormal basis for
the direct sum space, a {\em strong complement\/} to $\{x_j\}$. So
every Parseval frame has a strong complement which is unique up to
unitary equivalence. More generally, if $\{y_j\}$ is a general frame
we will call any frame $\{w_j\}$ such that $\{y_j\oplus w_j\}$ is a
Riesz basis for the direct sum space a {\em complementary\/} frame
(or {\em complement\/}) to $\{x_j\}$.

The notion of  strong complement has a natural generalization. Let
$\{x_n\}_{n\in \mathbb{J}}$ and $\{y_n\}_{n\in \mathbb{J}}$ be
Parseval frames in Hilbert spaces $H,K$, respectively, indexed by
the same set $\mathbb{J}$. Call these two frames {\em strongly
disjoint\/} if the (inner) direct sum $\{x_n\oplus y_n\colon \ n\in
\mathbb{J}\}$ is a Parseval frame for the direct sum Hilbert space
$H\oplus K$. It is not hard to see that this property of strong
disjointness is equivalent to the property that the ranges of their
analysis operators are orthogonal in $l^2(\mathbb{J})$. More
generally, we call a $k$-tuple of Parseval frames $(\{z_{1n}\}_{n\in
\mathbb{J}},\ldots, \{z_{kn}\}_{n\in\mathbb{J}})$ in Hilbert spaces
$H_1,\ldots, H_k$, respectively, a {\em strongly disjoint\/}
$k$-tuple if $\{z_{1n}\oplus \cdots \oplus z_{kn}\colon \ n\in
\mathbb{J}\}$ is a Parseval frame for $H_1\oplus \cdots \oplus H_k$,
and we call it a {\em complete\/} strongly disjoint $k$-tuple if
$\{z_{1n}\oplus\cdots\oplus z_{kn}\colon \ n\in \mathbb{J}\}$ is an
orthonormal basis for $H_1\oplus\cdots\oplus H_k$. If
$\theta_i\colon \ H_i\to l^2(\mathbb{J})$ is the frame transform,
$1\le i\le k$, then strong disjointness of a $k$-tuple is equivalent
to mutual orthogonality of $\{\text{ran } \theta_i\colon \ 1\le i\le
k\}$, and complete strong disjointness is equivalent to the
condition that $\bigoplus\limits^k_{i=1} \text{ran } \theta_i =
l^2(\mathbb{J})$.

There is a particularly simple intrinsic (i.e.\ non-geometric)
characterization of strong disjointness which is potentially useful
in applications:\ Let $\{x_n\}_{n\in \mathbb{J}}$ and
$\{y_n\}_{n\in\mathbb{J}}$ be Parseval frames for Hilbert spaces $H$
and $K$, respectively. Then $\{x_n\}$ and $\{y_n\}$ are strongly
disjoint if and only if one of the equations
\begin{align}\label{eq12}
&\sum_{n\in\mathbb{J}} \langle x,x_n\rangle y_n = 0\quad \text{for all } x\in
H\\
\text{or}\quad &\sum_{n\in\mathbb{J}} \langle y,y_n\rangle x_n= 0 \quad
\text{for all } y\in K\nonumber
\end{align}
holds. Moreover, if one holds the other holds also.

\subsection{Super-frames, Super-wavelets, and Multiplexing}
Suppose that $\{x_n\}_{n\in\mathbb{J}}$ and
$\{y_n\}_{n\in\mathbb{J}}$ are strongly disjoint Parseval frames for
Hilbert spaces $H$ and $K$, respectively. Then given any pair of
vectors $x\in H$, $y\in K$, we have that
\[
x = \sum_n \langle x,x_n\rangle x_n,\qquad y = \sum_n \langle y,y_n\rangle y_n.
\]
If  we let $a_n = \langle x,x_n\rangle$ and $b_n= \langle y,y_n\rangle$, and
then let $c_n=a_n+b_n$, we have
\[
\sum_n a_ny_n = 0,\qquad \sum_n b_nx_n = 0,
\]
by \eqref{eq12} and therefore we have
\begin{equation}\label{eq13}
x = \sum_n c_nx_n, \qquad y = \sum_N c_ny_n.
\end{equation}
This says that, by using one set of data $\{c_n\}$, we can recover
two vectors $x$ and $y$ (they may even lie in different Hilbert
spaces) by applying the respective inverse transforms (synthesis
operators) corresponding to the two frame $\{x_n\}$ and $\{y_n\}$.
The above argument obviously extends to the $k$-tuple case:\ If
$\{f_{in}\colon \ n\in \mathbb{J}\}$, $i=1,\ldots, k$, is a strongly
disjoint $k$-tuple of Parseval frames for Hilbert spaces
$H_1,\ldots, H_k$, and if $(x_1,\ldots, x_k)$ is an arbitrary
$k$-tuple of vectors with $x_i\in H_i$, $1\le i\le k$, then
\eqref{eq13} generalizes to
\[
x_i = \sum_{n\in \mathbb{J}} \langle x_i,f_{in}\rangle f_{in}
\]
for each $1\le i\le k$. So if we define a single ``master'' sequence of complex
numbers $\{c_n\colon \ n\in \mathbb{J}\}$ by
\[
c_n = \sum^k_{i=1} \langle x_i,f_{in}\rangle,
\]
then the strong disjointness implies that for each {\em individual\/} $i$ we
have
\[
x_i = \sum_{n\in \mathbb{J}} c_nf_{in}.
\]
This simple observation might be useful in applications to data
compression.

In [33] we called such an n-tuple of strongly disjoint (or simply
just disjoint) frames a \em{super-frame}, because it (or rather its
inner direct sum) is a frame for the \em{superspace} which s the
direct sum of the individual Hilbert spaces for the frames. In
connection with wavelet systems this observation lead us to the
notion of {\em superwavelet}, which is a particular type of
vector-valued wavelet. In operator-theoretic terms this is just a
restatement of the fact outlined above that a strongly disjoint
$k$-tuple of Parseval frames have frame-transforms which are
isometries into the same space $l^2(\mathbb{J})$ which have mutually
orthogonal ranges.

The notion of superframes and superwavelets, and many of their
properties, were also discovered and investigated by Radu Balan [3]
in his Ph.D. thesis, in work that was completely independent from
ours.

\subsection{Frame Vectors For Unitary Systems}
Let $\mathcal{U}$ be a unitary system on a Hilbert space $H$.
Suppose $\mathcal{W}(\mathcal{U})$ is nonempty, and fix $\psi\in
\mathcal{W}(\mathcal{U})$. Recall from Section~1 that if $\eta$ is
an arbitrary vector in $H $, then $\eta\in \mathcal{W}(\mathcal{U})$
if and only if there is a unitary $V$ (which is unique if it exists)
in the local commutant $\mathcal{C}_\psi(\mathcal{U})$ such that
$V\psi=\eta$. The following proposition shows that this idea
generalizes to the theory of frames. Analogously to the notion of a
wandering vector and a complete wandering vector, a vector $x\in H$
is called a {\em Parseval frame vector\/} (resp.\ {\em frame
vector\/} with bounds $a$ and $b$) for a unitary system
$\mathcal{U}$ if $\mathcal{U}x$ forms a tight frame (resp.\ frame
with bounds $a$ and $b$) for $\overline{span}(\mathcal{U}x)$. It is
called a {\em complete Parseval frame vector\/} (resp.\ {\em
complete frame vector\/} with bounds $a$ and $b$) when
$\mathcal{U}x$ is a Parseval frame (resp.\ frame with bounds $a$ and
$b$) for $H$.

\begin{proposition}\label{pro3.7}
Suppose that $\psi$ is a complete wandering vector for a unitary system
$\mathcal{U}$. Then
\begin{itemize}
\item[\em (i)] a vector $\eta$ is a Parseval frame vector for
$\mathcal{U}$ if and only if there is a (unique) partial isometry $A\in
C_\psi(\mathcal{U})$ such that $A\psi=\eta$.
\item[\em (ii)] a vector $\eta$ is a complete Parseval frame vector for
$\mathcal{U}$ if and only if there is a (unique) co-isometry $A\in
C_\psi(\mathcal{U})$ such that $A\psi=\eta$.
\end{itemize}
\end{proposition}

The above result does not tell the whole story. The reason is that
many unitary systems do not have wandering vectors but do have frame
vectors. For instance, this is the case in Example \ref{exmA}, where
the unitary system is the group of multiplication operators
$\mathcal{U} = \{M_{e^{ins}}\colon \ n \in \mathbb{Z}\}$ acting on
$L^2(E)$. In the case of a unitary system such as the wavelet system
$\mathcal{U}_{D,T}$ there exist {\em both\/} complete wandering
vectors {\em and\/} nontrivial Parseval frame vectors, so the theory
seems richer (however less tractable) and Proposition \ref{pro3.7}
is very relevant.

Much of Example \ref{exmA} generalizes to the case of an arbitrary countable
unitary group. There is a corresponding (geometric) dilation result.

\begin{proposition}\label{pro3.8}
Suppose that $\mathcal{U}$ is a unitary group such that
$\mathcal{W}(\mathcal{U})$ is non-empty. Then every complete
Parseval frame vector must be a complete wandering vector.
\end{proposition}

\begin{theorem}\label{thm3.9}
Suppose that $\mathcal{U}$ is a unitary group on $H$ and $\eta$ is a
complete Parseval frame vector for $\mathcal{U}$. Then there exists
a Hilbert space $K\supseteq H$ and a unitary group $\mathcal{G}$ on
$K$ such that $\mathcal{G}$ has complete wandering vectors, $H$ is
an invariant subspace of $\mathcal{G}$ such that $\mathcal{G}|_H =
\mathcal{U}$, and the map $g\to g|_H$ is a group isomorphism from
$\mathcal{G}$ onto $\mathcal{U}$.
\end{theorem}

The following is not hard, but it is very useful.

\begin{proposition}\label{pro3.10}
Suppose that $\mathcal{U}$ is a unitary group which has a complete
Parseval frame vector. Then the von Neumann algebra
$w^*(\mathcal{U})$ generated by $\mathcal{U}$ is finite.
\end{proposition}

\subsection{An Operator Model}
The following is a corollary of Theorem \ref{thm3.9}. It shows that Example
\ref{exmA} can be viewed as a model for certain operators.

\begin{corollary}\label{cor3.11}
Let $T\in B(H)$ be a unitary operator and let $\eta\in H$ be a
vector such that $\{T^n\eta\colon \ n\in \mathbb{Z}\}$ is a Parseval
frame for $H$. Then there is a unique (modulo a null set) measurable
set $E\subset  \mathbb{T}$ such that $\{T^n\eta\colon \ n\in
\mathbb{Z}\}$ and $\{e^{ins}|_E\colon \ n\in \mathbb{Z}\}$ are
unitarily equivalent frames.
\end{corollary}

\subsection{Group Representations}
These concepts generalize. For a unitary system $\mathcal{U}$ on a
Hilbert space $H$, a  closed subspace $M$ of $H$ is called a {\em
complete wandering subspace\/} for $\mathcal{U}$ if span $\{UM\colon
\ U\in \mathcal{U}\}$ is dense in $H$, and $UM\perp VM$ with $U\ne
V$. Let $\{e_i\colon \ i\in I\}$ be an orthonormal basis for $M$.
Then $M$ is a complete wandering subspace for $\mathcal{U}$ if and
only if $\{Ue_i\colon \ U\in \mathcal{U},i\in I\}$ is an orthonormal
basis for $H$. We call $\{e_i\}$ a {\em complete multi-wandering
vector}. Analogously, an $n$-tuple $(\eta_1,\ldots, \eta_n)$ of
non-zero vectors (here $n$ can be $\infty$) is called {\em complete
Parseval multi-frame vector\/} for $\mathcal{U}$ if $\{U\eta_i\colon
\ U\in \mathcal{U}$, $i=1,\ldots, n\}$ forms a complete Parseval
frame for $H$. Let $G$ be a group and let $\lambda$ be the left
regular representation of $G$ on $l^2(G)$. Then $\{\lambda_g\times
I_n\colon \ g\in G\}$ has a complete multi-wandering vector
$(f_1,\ldots, f_n)$, where $f_1 = (x_e,0,\ldots, 0),\ldots, f_n=$
$(0,0,\ldots, x_e)$. Let $P$ be any projection in the commutant of
$(\lambda\otimes I_n)(\mathcal{G})$. Then $(Pf_1,\ldots, Pf_n)$ is a
complete Parseval multi-frame vector for the subrepresentation
$(\lambda \otimes I_n)|_P$. It turns out that every representation
with a complete Parseval multi-frame vector arises in this way. Item
(i) of the following theorem is elementary and was mentioned
earlier; it is included for completeness.

\begin{theorem}\label{thm3.12}
Let $G$ be a countable group and let $\pi$ be a representation of $G$ on a
Hilbert space $H$. Let $\lambda$ denote the left regular representation of $G$
on $l^2(G)$. Then
\begin{itemize}
\item[\rm (i)] if $\pi(G)$ has a complete wandering vector then $\pi$ is
unitarily equivalent to $\lambda$,
\item[\rm (ii)] if $\pi(G)$ has a complete Parseval frame vector then
$\pi$ is unitarily equivalent to a subrepresentation of $\lambda$,
\item[\rm (iii)] if $\pi(G)$ has a complete Parseval multi-frame vector

$\{\psi_1,\psi_2,\ldots,\psi_n\}$, for some $1\le n<\infty$, then $\pi$ is
unitarily equivalent to a subrepresentation of $\lambda\otimes I_n$.
\end{itemize}
\end{theorem}

\section{Decompositions of Operators and Operator-Valued Frames}

The material we present here is contained in two recent papers.  The
first [15] was authored by a [VIGRE/REU] team consisting of K.
Dykema, D. Freeman, K. Kornelson, D. Larson, M. Ordower, and E.
Weber, with the
 title  \textit{Ellipsoidal Tight Frames}.
This article started as an undergraduate research project at Texas
A\&M in the summer of 2002, in which Dan Freeman was the student and
the other five were faculty mentors. Freeman is now a graduate
student at Texas A\&M. The project began as a solution of a finite
dimensional frame research problem, but developed into a rather
technically deep theory concerning a class of frames on an infinite
dimensional Hilbert space.  The second paper [44], entitled
\textit{Rank-one decomposition of operators and construction of
frames}, is a joint article by K. Kornelson and D. Larson.

\subsection{Ellipsoidal Frames}
We will use the term \textit{spherical frame} (or \emph{uniform
frame}) for a frame sequence which is \textit{uniform} in the sense
that all its vectors have the same norm.  Spherical frames which are
tight have been the focus of several articles by different
researchers.  Since frame theory is essentially geometric in nature,
from a purely mathematical point of view it is natural  to ask:
Which other surfaces in a finite or infinite dimensional Hilbert
space contain tight frames?  (These problems can make darn good REU
projects, in particular.) In the first article we considered
ellipsoidal surfaces.

By an \textit{ellipsoidal surface} we mean the image of the unit
sphere $S_1$ in the underlying Hilbert space $H$ under a bounded
invertible operator $A$ in $B(H)$, the set of all bounded linear
operators on $H$.  Let $E_A$ denote the ellipsoidal surface $E_A :=
AS_1$.  A frame contained in $E_A$ is called an  \textit{ellipsoidal
frame}, and if it is tight it is called an ellipsoidal tight frame
(ETF) for that surface.  We say that a frame bound $K$ is
\textit{attainable} for $E_A$ if there is an ETF for $E_A$ with
frame bound K.

Given an ellipsoidal surface $E := E_A$, we can assume $E = E_T$
where T is a positive invertible operator.  Indeed, given an
invertible operator $A$, let $A^* = U|A^*|$ be the polar
decomposition, where $|A^*| = (AA^*)^{1/2}$.  Then $A = |A^*|U^*$.
By taking $T = |A^*|$, we see that $TS_1 = AS_1$.  Moreover, it is
easily seen that the positive operator $T$ for which $E = E_T$ is
unique.

The starting point for the work in the first paper was the following
Proposition.  For his REU project Freeman found an elementary
calculus proof of this for the real case.  Others have also
independently found this result, including V. Paulsen, and P.
Casazza and M. Leon.

\begin{proposition} Let $E_A$ be an ellipsoidal surface on a finite dimensional
real or complex
Hilbert space $H$ of dimension $n$.   Then for any integer $k \geq
n$, $E_A$ contains a tight frame of length $k$, and every ETF on
$E_A$ of length $k$ has frame bound $K = k
\left[\text{trace}(T^{-2})\right]^{-1}$.
 \end{proposition}
We use the following standard definition: For an operator $B \in H$,
the \textit{essential norm} of $B$ is:
$$ \|B\|_{ess} := \inf \{\|B-K\|\;:\;  \text{$K$ is a compact operator in}
B(H)\} $$
Our main frame theorem from the first paper is:

\begin{theorem} Let $E_A$ be an ellipsoidal surface in an infinite dimensional
real or complex
Hilbert space.  Then for any constant
 $K >   \|T^{-2}\|^{-1}_{ess} $,
 $E_T$ contains a tight frame with frame bound $K$.
 \end{theorem}

So, for fixed $A$, in finite dimensions the set of attainable ETF
frame bounds is finite, whereas in infinite dimensions it is a
continuum.

\textit{Problem}.  If the essential norm of $A$ is replaced with the
norm of $A$ in the above theorem, or if the inequality is replaced
with equality, then except for some special cases, and trivial
cases, no theorems of any degree of generality are known concerning
the set of attainable frame bounds for ETF's on $E_A$.  It would be
interesting to have a general analysis of the case where $A - I$ is
compact.  In this case, one would want to know necessary and
sufficient conditions for existence of a tight frame on $E_A$ with
frame bound 1. In the special case $A = I$ then, of course, any
orthonormal basis will do, and these are the only tight frames on
$E_A$ in this case.  What happens in general when $ \|A\|_{ess} = 1
$ and $A$ is a small perturbation of $I$?

We use elementary tensor notation for a rank-one operator on $H$.
Given $u,v,x \in H$, the operator $u \otimes v$ is defined by $(u
\otimes v)x = \langle x,v \rangle u$ for $x \in H$. The operator $u
\otimes u$ is a projection if and only if $\|u\|=1$.

Let $\{x_j\}_j $ be a frame for $H$.
 The standard
frame operator is defined by: $ Sw = \sum_j \langle w,x_j \rangle
x_j  = \sum_j
  \left( x_j \otimes x_j \right)w $ .
\noindent Thus $S = \sum_j x_j \otimes x_j$, where this series of
positive rank-1
 operators converges in the
strong operator topology (i.e. the topology of pointwise
convergence). In the
 special case where each $x_j$ is a unit vector, $S$ is the sum of the rank-1
projections $P_j = x_j \otimes x_j$.

For $A$ a positive operator, we say that $A$ has a {\it projection
decomposition} if $A$ can be expressed as the sum of a finite or
infinite sequence of (not necessarily mutually orthogonal)
self-adjoint projections, with convergence in the strong operator
topology.

If  $x_j$ is a frame of unit vectors, then $S = \sum_j x_j \otimes
x_j$ is a projection decomposition of the frame operator. This
argument is trivially reversible, so a positive invertible operator
$S$ is the
 frame operator for a frame
of unit vectors if and only if it admits a projection decomposition
$S = \sum_j P_J$. If the projections in the decomposition are not of
rank one, each projection can be further decomposed (orthogonally)
 into rank-1
projections, as needed, expressing $S = \sum_n x_n \otimes x_n$, and
then the sequence $\{x_n\}$ is a frame of unit vectors with frame
operator $S$.

In order to prove Theorem 22, we first proved Theorem 23 (below),
using purely operator-theoretic techniques.

 \begin{theorem}  Let $A$ be a positive operator in $B(H)$ for $H$ a real or
complex Hilbert space with infinite dimension, and suppose
 $\|A\|_{ess}>1$.
Then $A$ has a projection decomposition.
\end{theorem}
Suppose, then, that $\{x_n\}$ is a frame of unit vectors with frame
operator $S$. If we let $y_j = S^{-\frac 12}x_j$, then
 $\{y_j\}_j $ is a \textit{Parseval} frame.
 So  $\{y_j\}_j$ is an
ellipsoidal tight frame for the ellipsoidal surface $E_{S^{-\frac
12}}=S^{-\frac 12}S_1$.  This argument is reversible: Given a
positive invertible operator $T$, let $S = T^{-2}$.  Scale $T$ if
necessary so that $\|S\|_{ess}>1$.  Let $S = \sum_j x_j \otimes x_j$
be a projection decomposition of $S$. Then $\{Tx_j\}$ is an ETF for
the ellipsoidal surface $TS_1$.  Consideration of frame bounds and
scale factors then yields Theorem 22.

Most of our second paper concerned \textit{weighted} projection
decompositions of positive operators, and resultant theorems
concerning frames.  If $T$ is a positive operator, and if $\{c_n\}$
is a sequence of positive scalars, then a weighted projection
decomposition of $T$ with weights $\{c_n\}$ is a decomposition $T =
\sum_j P_j$ where the $P_j$ are projections, and the series
converges strongly. We have since adopted the term \textit{targeted}
to refer to such a decomposition, and generalizations thereof.  By a
\textit{targeted decomposition} of $T$ we mean any strongly
convergent decomposition  $T = \sum_n T_n$ where the $T_n$ is a
sequence of \textit{simpler} positive operators with special
prescribed properties. So a weighted decomposition is a targeted
decomposition for which the scalar weights are the prescribed
properties.  And, of course, a projection decomposition is a special
case of targeted decomposition.

After a sequence of Lemmas, building up from finite dimensions and
employing spectral theory for operators, we arrived at the following
theorem.  We will not discuss the details here because of limited
space. It is the \textit{weighted} analogue of theorem 23.

\begin{theorem}  Let $B$ be a positive operator in $B(H)$ for $H$ with
$\|B\|_{ess}>1$.
Let $\{c_i\}_{i=1}^{\infty}$ be any sequence of numbers with $0<c_i
\leq 1$ such that $\sum_i c_i = \infty$. Then there exists a
sequence of rank-one projections $\{P_i\}_{i=1}^{\infty}$ such that
$ B = \sum_{i=1}^{\infty} c_i P_i $
\end{theorem}

\bigskip


\subsection{A Problem in Operator Theory}
We will discuss a problem in operator theory that was motivated by a
problem in the theory of Modulation Spaces.  We tried to obtain an
actual "reformulation" of the modulation space problem in terms of
operator theory,  and it is well possible that such a reformulation
can be found.  At the least we (Chris Heil and myself) found the
following operator theory problem, whose solution could conceivably
impact mathematics beyond operator theory.  I find it rather
fascinating.  I need to note that we subsequently showed (in an
unpublished jointly-written expository article) that the actual
\emph{modulation space} connection requires a modified and more
sophisticated version of the problem we present below.  I still
feel, that the problem I will present here has some independent
interest, and may serve as a "first step" in developing a theory
that might have some usefulness.  Thus, I hope that the reader will
find it interesting.

\bigskip
\bigskip
\noindent

Let $H$ be an infinite dimensional separable Hilbert space. As
usual, denote the Hilbert space norm on H by $\| \cdot \|$. If $x$
and $y$ are vectors in $H$, then  $x \otimes  y$ will denote the
operator of rank one defined by $(x \otimes y)z =  \langle z, y
\rangle x$. The operator norm of $x \otimes y$ is then just the
product of $\|x\|$ and $\|y\|$.

\medskip

Fix an orthonormal basis $\{e_n\}_n$ for $H$. For each vector $v$ in
$H$, define  $$\||v\|| = \sum_n |\langle v,e_n \rangle |$$ This may
be $+\infty$.

\medskip

Let  $L$  be the set of all vectors $v$ in H for which  $\||v\||$ is
finite.  Then $L$ is a dense linear subspace of H, and is a Banach
space in the "triple norm".  It is of    course isomorphic to
$\ell^1$

\medskip

Let $T$ be any positive trace-class operator in $B(H)$.

The usual eigenvector decomposition  for $T$ expresses  $T$  as a
series converging in the strong operator topology of operators $h_n
\otimes h_n$,  where $\{h_n\}$ is an orthogonal sequence of
eigenvectors of $T$.  That is,
$$T = \sum_n h_n \otimes h_n$$  In this representation the eigenvalue
corresponding to the eigenvector $h_n$ is the square of the norm: $\|h_n\|^2$.
The trace of  $T$  is then  $$\sum_n \|h_n\|^2$$  and since $T$ is positive this
is also the trace-class norm of $T$.

\medskip

Let us say that  $T$  is of  \emph{Type A}  with respect to the
orthonormal basis $\{e_n\}$ if, for  the eigenvectors $\{h_n\}$  as
above,  we have  that $\sum_n \||h_n\||^2$ is finite.  [Note that
this is just the (somewhat unusual) formula displayed above for the
trace of T with the triple norm used in place of the usual Hilbert
space norm of the vectors  $\{h_n\}$.]

\medskip

And let us say that  $T$  is of  \emph{Type B}  with respect to the
orthonormal basis $\{e_n\}$ if there is \emph{some} sequence of
vectors $\{v_n\}$ in H with $\sum_n \||v_n\||^2$  \emph{finite}
such that
$$T = \sum_n v_n \otimes h_n$$ where the convergence of this series
is in the strong operator topology.

\bigskip

\textbf{Problem:}   If  $T$  if of \emph{Type B} with respect to an
orthonormal basis $\{e_n\}$, then must it be of \emph{Type A} with
respect to $\{e_n\}$?

\medskip

\textit{Note:} If the answer to this problem is negative (as I
suspect it is), then the following subproblem would be an
interesting one.

\bigskip

\textbf{Subproblem:}  Let $\{e_n\}$ be an orthonormal basis for $H$.
Find a characterization of all positive trace class operators  $T$
that are of \emph{Type B}  with respect to $\{e_n\}$.  In
particular, is every positive trace class operator $T$ of \emph{Type
B} with respect to $\{e_n\}$?  My feeling is \emph{no}. (See the
next example.)

\medskip

\begin{example}\label{ex2}  Let  $x$ be any vector in $H$ that is not in $L$,
and let $T = x
\otimes x$.  Then $T$ is trace class, in fact has rank one, but
clearly $T$ is clearly \emph{not} of Type A.  Can such a $T$ be of
type B?  (I don't think it is necessarily of Type B for all such
$T$, however.)
\end{example}

\end{document}